\documentclass[a4paper,11pt]{article}
\usepackage{amsfonts,amssymb,amsmath}

\usepackage[english]{babel}
 \makeatletter
 \@addtoreset{equation}{section}
 \makeatother

 \makeatletter
 \@addtoreset{enunciato}{section}
 \makeatother

 \newcounter{enunciato}[section]

 \newtheorem{ittheorem}{Theorem}
 \newtheorem{itlemma}{Lemma}
 \newtheorem{itproposition}{Proposition}
 \newtheorem{itdefinition}{Definition}
 \newtheorem{itremark}{Remark}
 \newtheorem{itclaim}{Claim}
 \newtheorem{itfact}{Fact}
 \newtheorem{itconjecture}{Conjecture}

 \newenvironment{theorem}{\addtocounter{enunciato}{1}
 \begin{ittheorem}}{\end{ittheorem}}

 \newenvironment{lemma}{\addtocounter{enunciato}{1}
 \begin{itlemma}}{\end{itlemma}}

 \newenvironment{proposition}{\addtocounter{enunciato}{1}
 \begin{itproposition}}{\end{itproposition}}

 \newenvironment{definition}{\addtocounter{enunciato}{1}
 \begin{itdefinition}}{\end{itdefinition}}

 \newenvironment{remark}{\addtocounter{enunciato}{1}
 \begin{itremark}}{\end{itremark}}

 \newenvironment{claim}{\addtocounter{enunciato}{1}
 \begin{itclaim}}{\end{itclaim}}

 \newenvironment{fact}{\addtocounter{enunciato}{1}
 \begin{itfact}}{\end{itfact}}

 \newenvironment{conjecture}{\addtocounter{enunciato}{1}
 \begin{itconjecture}}{\end{itconjecture}}

 \newcommand{\be}[1]{\begin{equation}\label{#1}}
 \newcommand{\ee}{\end{equation}}

 \newcommand{\bl}[1]{\begin{lemma}\label{#1}}
 \newcommand{\el}{\end{lemma}}

 \newcommand{\br}[1]{\begin{remark}\label{#1}}
 \newcommand{\er}{\end{remark}}

 \newcommand{\bt}[1]{\begin{theorem}\label{#1}}
 \newcommand{\et}{\end{theorem}}

 \newcommand{\bd}[1]{\begin{definition}\label{#1}}
 \newcommand{\ed}{\end{definition}}

 \newcommand{\bcl}[1]{\begin{claim}\label{#1}}
 \newcommand{\ecl}{\end{claim}}

 \newcommand{\bfact}[1]{\begin{fact}\label{#1}}
 \newcommand{\efact}{\end{fact}}

 \newcommand{\bp}[1]{\begin{proposition}\label{#1}}
 \newcommand{\ep}{\end{proposition}}

 \newcommand{\bc}[1]{\begin{corollary}\label{#1}}
 \newcommand{\ec}{\end{corollary}}

 \newcommand{\bcj}[1]{\begin{conjecture}\label{#1}}
 \newcommand{\ecj}{\end{conjecture}}

 \newcommand{\bpr}{\begin{proof}}
 \newcommand{\epr}{\end{proof}}

 \newcommand{\bprl}[1]{\begin{proofof}{\it\ref{#1}}.\,\,}
 \newcommand{\eprl}{\end{proofof}}

 \newcommand{\bi}{\begin{itemize}}
 \newcommand{\ei}{\end{itemize}}

 \newcommand{\ben}{\begin{enumerate}}
 \newcommand{\een}{\end{enumerate}}

 \newenvironment{proof}{\noindent {\em Proof}.\,\,}{\hspace*{\fill}$\halmos$\medskip}
 \newenvironment{proofof}{\noindent {\em Proof of Lemma\,\,}}{\hspace*{\fill}$\halmos$\medskip}
 \newcommand{\halmos}{\rule{1ex}{1.4ex}}

 \newcommand{\one}{{\mathchoice {1\mskip-4mu\mathrm l}
         {1\mskip-4mu\mathrm l}
         {1\mskip-4.5mu\mathrm l}
         {1\mskip-5mu\mathrm l}}}

\def \N {{\mathbb N}}

\def \R {{\mathbb R}}
\def \Z {{\mathbb Z}}

\def \ra {\rightarrow}
\def \ba {\begin{array}}
\def \ea {\end{array}}
\def \da {\downarrow}
\def \ua {\uparrow}

\def \P {{\mathbb P}}
\def \E {{\mathbb E}}

\def \ES {{\rm E}}
\def \PS {{\rm P}}
\def \Sp {{\rm Sp}}

\def \new {{\rm new}}

\def \c {{\rm c}}

\def \per {{\rm per}}
\def \onek {{1[\kappa]}}

\def \cA {{\mathcal A}}

\def \cC {{\mathcal C}}
\def \cD {{\mathcal D}}
\def \cE {{\mathcal E}}

\def \cG {{\mathcal G}}

\def \cK {{\mathcal K}}

\def \cM {{\mathcal M}}
\def \cP {{\mathcal P}}
\def \cQ {{\mathcal Q}}

\def \cS {{\mathcal S}}
\def \cT {{\mathcal T}}

\def\one{\rlap{\mbox{\small\rm 1}}\kern.15em 1}

\begin{document}
\title{Intermittency on catalysts:\\
three-dimensional simple symmetric exclusion}

\author{\renewcommand{\thefootnote}{\arabic{footnote}}
J.\ G\"artner
\footnotemark[1]
\\
\renewcommand{\thefootnote}{2,3}
F.\ den Hollander
\footnotemark
\\
\renewcommand{\thefootnote}{\arabic{footnote}}
G.\ Maillard
\footnotemark[4]
}

\footnotetext[1]{
Institut f\"ur Mathematik, Technische Universit\"at Berlin,
Strasse des 17.\ Juni 136, D-10623 Berlin, Germany,
{\sl jg@math.tu-berlin.de}
}
\footnotetext[2]{
Mathematical Institute, Leiden University, P.O.\ Box 9512,
2300 RA Leiden, The Netherlands,\\
{\sl denholla@math.leidenuniv.nl}
}
\footnotetext[3]{
EURANDOM, P.O.\ Box 513, 5600 MB Eindhoven, The Netherlands
}
\footnotetext[4]{
CMI-LATP, Universit\'e de Provence,
39 rue F. Joliot-Curie, F-13453 Marseille Cedex 13, France,\\
{\sl maillard@cmi.univ-mrs.fr}
}

\date{16th December 2008}
\maketitle

\begin{abstract}
We continue our study of intermittency for the parabolic Anderson model
$\partial u/\partial t = \kappa\Delta u + \xi u$ in a space-time random
medium $\xi$, where $\kappa$ is a positive diffusion constant, $\Delta$
is the lattice Laplacian on $\Z^d$, $d \geq 1$, and $\xi$ is a simple
symmetric exclusion process on $\Z^d$ in Bernoulli equilibrium. This model
describes the evolution of a \emph{reactant} $u$ under the influence of
a \emph{catalyst} $\xi$.

In \cite{garholmai06} we investigated the behavior of the annealed Lyapunov
exponents, i.e., the exponential growth rates as $t\to\infty$ of the successive
moments of the solution $u$. This led to an almost complete picture of intermittency
as a function of $d$ and $\kappa$. In the present paper we finish our study by
focussing on the asymptotics of the Lyaponov exponents as $\kappa\to\infty$ in
the \emph{critical} dimension $d=3$, which was left open in \cite{garholmai06}
and which is the most challenging. We show that, interestingly, this asymptotics
is characterized not only by a \emph{Green} term, as in $d\geq 4$, but also by a
\emph{polaron} term. The presence of the latter implies intermittency of \emph{all}
orders above a finite threshold for $\kappa$.

{\it Acknowledgment.} The research of this paper was partially supported by the
DFG Research Group 718 ``Analysis and Stochastics in Complex Physical Systems'',
the DFG-NWO Bilateral Research Group ``Mathematical Models from Physics and
Biology'', and the ANR-project MEMEMO.
\end{abstract}

\newpage


\section{Introduction and main result}
\label{S1}


\subsection{Model}
\label{S1.1}

In this paper we consider the \emph{parabolic Anderson model} (PAM) on $\Z^d$, $d \geq 1$,
\be{pA}
\begin{cases}
\displaystyle\frac{\partial u}{\partial t} = \kappa\Delta u + \xi u
&\text{on }\Z^d\times[0,\infty),\\
u(\cdot,0) = 1
&\text{on }\Z^d,
\end{cases}
\ee
where $\kappa$ is a positive diffusion constant, $\Delta$ is the lattice Laplacian
acting on $u$ as
\be{dL}
\Delta u(x,t) = \sum_{{y\in\Z^d} \atop {\|y-x\|=1}} [u(y,t)-u(x,t)]
\ee
($\|\cdot\|$ is the Euclidian norm), and
\be{xigen}
\xi=(\xi_t)_{t \geq 0}, \qquad \xi_t = \{\xi_t(x)\colon\,x\in\Z^d\},
\ee
is a space-time random field that drives the evolution. If $\xi$ is given by an infinite
particle system dynamics, then the solution $u$ of the PAM may be interpreted as the
concentration of a diffusing \emph{reactant} under the influence of a \emph{catalyst}
performing such a dynamics.

In G\"artner, den Hollander and Maillard~\cite{garholmai06} we studied the PAM for $\xi$
Symmetric Exclusion (SE), and developed an almost complete qualitative picture. In the
present paper we finish our study by focussing on the limiting behavior as $\kappa\to\infty$
in the \emph{critical} dimension $d=3$, which was left open in \cite{garholmai06} and which
is the most challenging. We restrict to \emph{Simple Symmetric Exclusion} (SSE), i.e.,
$(\xi_t)_{t\geq 0}$ is the Markov dynamics on $\Omega=\{0,1\}^{\Z^3}$ ($0$ = vacancy,
$1$ = particle) with generator $L$ acting on cylinder functions $f\colon\Omega\to\R$ as
\be{expro1}
(Lf)(\eta)
= \frac16 \sum_{\{a,b\}}
\Big[f\big(\eta^{a,b}\big)-f(\eta)\Big],
\quad\eta\in\Omega,
\ee
where the sum is taken over all unoriented nearest-neighbor bonds $\{a,b\}$ of $\Z^3$, and
$\eta^{a,b}$ denotes the configuration obtained from $\eta$ by interchanging the states at
$a$ and $b$:
\be{SSEflip}
\eta^{a,b}(a)=\eta(b), \quad \eta^{a,b}(b)=\eta(a), \quad  \eta^{a,b}(x)=\eta(x)
\mbox{ for } x\notin\{a,b\}.
\ee
(See Liggett~\cite{lig85}, Chapter VIII.) Let $\P_\eta$ and $\E_\eta$ denote probability and
expectation for $\xi$ given $\xi_0=\eta\in\Omega$. Let $\xi_0$ be drawn according to the
Bernoulli product measure $\nu_\rho$ on $\Omega$ with density $\rho\in(0,1)$. The probability
measures $\nu_\rho$, $\rho\in(0,1)$, are the only extremal equilibria of the SSE dynamics.
(See Liggett~\cite{lig85}, Chapter VIII, Theorem 1.44.) We write $\P_{\nu_\rho} = \int_\Omega
\nu_{\rho}(d\eta)\,\P_{\eta}$ and $\E_{\,\nu_\rho} = \int_\Omega \nu_{\rho}(d\eta)\,\E_{\,\eta}$.


\subsection{Lyapunov exponents}
\label{S1.2}

For $p\in\N$, define the $p$-th \emph{annealed Lyapunov exponent} of the PAM by
\be{lyapdef}
\lambda_p(\kappa,\rho)
= \lim_{t\to\infty}\frac{1}{pt} \log \E_{\,\nu_\rho}\left([u(0,t)]^p\right).
\ee
We are interested in the asymptotic behavior of $\lambda_p(\kappa,\rho)$ as $\kappa\to\infty$
for fixed $\rho$ and $p$. To this end, let $G$ denote the value at $0$ of the \emph{Green
function} of simple random walk on $\Z^3$ with jump rate $1$ (i.e., the Markov process with
generator $\frac16\Delta$), and let $\cP_3$ be the value of the \emph{polaron variational problem}
\be{Rpdef}
\cP_3 = \sup_{{f \in H^1(\R^3)} \atop {\|f\|_2=1}}
\left[\,\left\|\left(-\Delta_{\R^3}\right)^{-1/2}\,f^2\right\|_2^2
- \left\|\nabla_{\R^3} f\right\|_2^2\,\right],
\ee
where $\nabla_{\R^3}$ and $\Delta_{\R^3}$ are the \emph{continuous} gradient and
Laplacian, $\|\cdot\|_2$ is the $L^2(\R^3)$-norm, $H^1(\R^3)=\{f\in L^2(\R^3)\colon\,
\nabla_{\R^3}f\in L^2(\R^3)\}$, and
\be{Gint}
\left\|\left(-\Delta_{\R^3}\right)^{-1/2}\,f^2\right\|_2^2
= \int_{\R^3} dx\,f^2(x) \int_{\R^3} dy\,f^2(y)\,\,\frac{1}{4\pi\|x-y\|}.
\ee
(See Donsker and Varadhan~\cite{donvar83} for background on how $\cP_3$ arises in
the context of a self-attracting Brownian motion referred to as the polaron model.
See also G\"artner and den Hollander~\cite{garhol04}, Section 1.5.)

We are now ready to formulate our main result (which was already announced in
G\"artner, den Hollander and Maillard~\cite{garholmaiHvW}).

\bt{main}
Let $d=3$, $\rho\in (0,1)$ and $p\in\N$. Then
\be{limlamb*}
\lim_{\kappa\to\infty} \kappa[\lambda_p(\kappa,\rho)-\rho] =
\frac{1}{6}\, \rho(1-\rho)G + [6\rho(1-\rho)p]^2\cP_3.
\ee
\et

Note that the expression in the r.h.s.\ of (\ref{limlamb*}) is the sum of a \emph{Green} term
and a \emph{polaron} term. The existence, continuity, monotonicity and convexity of $\kappa
\mapsto\lambda_p(\kappa,\rho)$ were proved in \cite{garholmai06} for all $d \geq 1$ for all
exclusion processes with an irreducible and symmetric random walk transition kernel. It was
further proved that $\lambda_p(\kappa,\rho)=1$ when the random walk is recurrent and
$\rho<\lambda_p(\kappa,\rho)<1$ when the random walk is transient. Moreover, it was shown that
for simple random walk in $d\geq 4$ the asymptotics as $\kappa\to\infty$ of $\lambda_p(\kappa,\rho)$
is similar to (\ref{limlamb*}), but \emph{without} the polaron term. In fact, the subtlety in $d=3$
is caused by the appearance of this extra term which, as we will see in Section~\ref{S5}, is related
to the large deviation behavior of the occupation time measure of a rescaled random walk that lies
deeply hidden in the problem. For the heuristics behind Theorem \ref{main} we refer the reader to
\cite{garholmai06}, Section 1.5.


\subsection{Intermittency}
\label{S1.3}

The presence of the polaron term in Theorem~\ref{main} implies that, for each $\rho\in(0,1)$,
there exists a $\kappa_0(\rho)>0$ such that the \emph{strict} inequality
\be{intermiteq}
\lambda_p(\kappa,\rho)>\lambda_{p-1}(\kappa,\rho)
\quad\forall\,\kappa>\kappa_0(\rho)
\ee
holds for $p=2$ and, consequently, for all $p \geq 2$ by the convexity of $p\mapsto
p\,\lambda_p(\kappa,\rho)$. This means that all moments of the solution $u$ are \emph{intermittent}
for $\kappa>\kappa_0(\rho)$, i.e., for large $t$  the random  field $u(\cdot,t)$ develops
sparse high spatial peaks dominating the moments in such a way that each moment is dominated
by its own collection of peaks (see G\"artner and K\"onig  \cite{garkon04}, Section 1.3, and
den Hollander \cite{hol00}, Chapter 8, for more explanation).

In \cite{garholmai06} it was shown that for all $d\geq 3$ the PAM is intermittent for \emph{small}
$\kappa$. We conjecture that in $d=3$ it is in fact intermittent for \emph{all} $\kappa$.
Unfortunately, our analysis does not allow us to treat intermediate values of $\kappa$ (see
the figure).


\setlength{\unitlength}{0.25cm}
\begin{picture}(20,14)(-15,4)
\put(0,-2){\line(18,0){18}}
\put(0,-2){\line(0,15){15}}
{\thicklines
\qbezier(0,8)(2,6.6)(4,5.3)
\qbezier(0,6)(2,4.8)(4,3.8)
\qbezier(0,4)(2,3)(4,2.25)
\qbezier(9,2.9)(12,1.3)(15,0.9)
\qbezier(9,2)(12,1)(15,0.6)
\qbezier(9,1.1)(12,0.7)(15,0.3)
}
\qbezier[60](0,0)(9,0)(17,0)
\qbezier[60](0,10.5)(9,10.5)(17,10.5)
\put(-1,-3.5){$0$}
\put(-1.4,-0.3){$\rho$}
\put(-1.4,10.2){$1$}
\put(0,8){\circle*{.45}}
\put(0,6){\circle*{.45}}
\put(0,4){\circle*{.45}}
\put(-4.3,8){\small $p=3$}
\put(-4.3,6){\small $p=2$}
\put(-4.3,4){\small $p=1$}
\put(6,2.5){{\bf ?}}
\put(19,-2.3){$\kappa$}
\put(-1.5,15){$\lambda_p(\kappa)$}
\put(-4,-7){\small
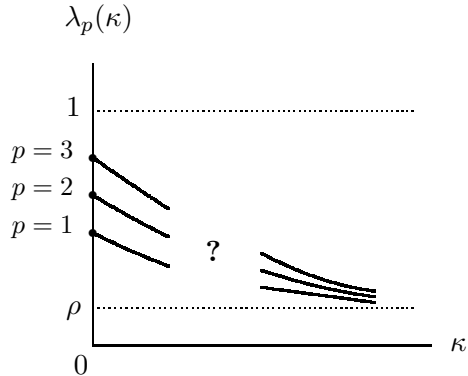
Qualitative picture of $\kappa\mapsto\lambda_p(\kappa)$ for $p=1,2,3$.}
\end{picture}


\vspace{3cm}

The formulation of Theorem~\ref{main} coincides with the corresponding result in G\"artner
and den Hollander~\cite{garhol04}, where the random potential $\xi$ is given by independent
simple random walks in a Poisson equilibrium in the so-called weakly catalytic regime. However,
as we already pointed out in \cite{garholmai06}, the approach in \cite{garhol04} cannot be
adapted to the exclusion process, since it relies on an explicit Feynman-Kac representation
for the moments that is available only in the case of \emph{independent} particle motion. We
must therefore proceed in a totally different way. Only at the end of Section~\ref{S5} will
we be able to use some of the ideas in \cite{garhol04}.


\subsection{Outline}
\label{S1.4}

Each of Sections~\ref{S2}--\ref{S5} is devoted to a major step in the proof of
Theorem~\ref{main} for $p=1$. The extension to $p \geq 2$ will be indicated in
Section~\ref{S6}.

In Section \ref{S2} we start with the Feynman-Kac representation for the first moment
of the solution $u$, which involves a random walk sampling the exclusion process. After
rescaling time, we transform the representation w.r.t.\ the old measure to a representation
w.r.t.\ a new measure via an appropriate absolutely continuous transformation. This allows
us to separate the parts responsible for, respectively, the Green term and the polaron term
in the r.h.s.\ of (\ref{limlamb*}). Since the Green term has already been handled in
\cite{garholmai06}, we need only concentrate on the polaron term. In Section~\ref{S3} we
show that, in the limit as $\kappa\to\infty$, the new measure may be replaced by the old
measure. The resulting representation is used in Section \ref{S4} to prove the \emph{lower
bound} for the polaron term. This is done analytically with the help of a Rayleigh-Ritz
formula. In Section~\ref{S5}, which is technical and takes up almost half of the paper, we
prove the corresponding \emph{upper bound}. This is done by freezing and defreezing
the exclusion process over long time intervals, allowing us to approximate the representation
in terms of the occupation time measures of the random walk over these time intervals. After
applying spectral estimates and using a large deviation principle for these occupation time
measures, we arrive at the polaron variational formula.


\section{Separation of the Green term and the polaron term}
\label{S2}

In Section~\ref{S2.1} we formulate the Feynman-Kac representation for the first moment of
$u$ and show how to split this into two parts after an appropriate change of measure.
In Section~\ref{S2.2} we formulate two propositions for the asymptotics of these two parts,
which lead to, respectively, the Green term and the polaron term in (\ref{limlamb*}). These
two propositions will be proved in Sections~\ref{S3}--\ref{S5}. In Section~\ref{S2.3} we
state and prove three elementary lemmas that will be needed along the way.


\subsection{Key objects}
\label{S2.1}

The solution $u$ of the PAM in (\ref{pA}) admits the Feynman-Kac representation
\be{fey-kac1}
u(x,t) = \ES_{\,x}^X\left(\exp\left[\int_0^t ds\,\,
\xi_{t-s}\left(X_{\kappa s}\right)\right]\right),
\ee
where $X$ is simple random walk on $\Z^3$ with step rate $6$ (i.e., with generator $\Delta$)
and $\PS_{x}^X$ and $\ES_{\,x}^X$ denote probability and expectation with respect to $X$
given $X_0=x$. Since $\xi$ is reversible w.r.t\ $\nu_\rho$, we may reverse time in (\ref{fey-kac1})
to obtain
\be{expectu}
\E_{\,\nu_\rho}\big(u(0,t)\big)
=\E_{\,\nu_\rho,0}\bigg(\exp\bigg[\int_0^t ds\,\, \xi_s\big(X_{\kappa s}\big)\bigg]\bigg),
\ee
where $\E_{\,\nu_\rho,0}$ is expectation w.r.t.\ $\P_{\nu_\rho,0}=\P_{\nu_\rho}\otimes \PS_0^X$.

As in \cite{garhol04} and \cite{garholmai06}, we rescale time and write
\be{expectu*}
e^{-\rho(t/\kappa)}\E_{\,\nu_\rho}\big(u(0,t/\kappa)\big)
=\E_{\,\nu_\rho,0}\bigg(\exp\bigg[\frac1\kappa\int_0^t ds\, \phi(Z_s)\bigg]\bigg)
\ee
with
\be{phidef}
\phi(\eta,x)=\eta(x)-\rho
\ee
and
\be{Ztdef}
Z_t = \big(\xi_{t/\kappa},X_t\big).
\ee
From (\ref{expectu*}) it is obvious that (\ref{limlamb*}) in Theorem \ref{main} (for $p=1$)
reduces to
\be{limscalfin}
\lim_{\kappa\to\infty} \kappa^2\lambda^*(\kappa)
= \frac1{6}\, \rho(1-\rho)G + [6\rho(1-\rho)]^2\cP_3,
\ee
where
\be{lscal}
\lambda^*(\kappa) = \lim_{t\to\infty} \frac{1}{t} \log
\E_{\,\nu_\rho,0}\left(\exp\left[\frac{1}{\kappa}\int_0^t ds\,
\phi(Z_s)\right]\right).
\ee
Here and in the rest of the paper we suppress the dependence on $\rho\in(0,1)$ from
the notation. Under $\P_{\eta,x}=\P_{\eta}\otimes\PS_{x}^X$, $(Z_t)_{t\geq 0}$ is
a Markov process with state space $\Omega\times\Z^3$ and generator
\be{lb2.7}
\cA = \frac{1}{\kappa} L + \Delta
\ee
(acting on the Banach space of bounded continuous functions on $\Omega\times \Z^3$,
equipped with the supremum norm). Let $(\cS_t)_{t\geq 0}$ denote the semigroup generated
by $\cA$.

Our aim is to make an absolutely continuous transformation of the measure $\P_{\eta,x}$
with the help of an exponential martingale, in such a way that, under the new measure
$\P_{\eta,x}^\new$, $(Z_t)_{t\geq 0}$ is a Markov process with generator $\cA^\new$
of the form
\be{lb2i.2}
\cA^{\new}f
= e^{-\frac{1}{\kappa}\psi}\cA\left(e^{\frac{1}{\kappa}\psi}f
\right) - \left(e^{-\frac{1}{\kappa}\psi}\cA e^{\frac{1}{\kappa}\psi}
\right)f.
\ee
This transformation leads to an interaction between the exclusion process part and the
random walk part of $(Z_t)_{t\geq 0}$, controlled by $\psi\colon \Omega\times\Z^3\to\R$.
As explained in \cite{garholmai06}, Section 4.2, it will be expedient to choose $\psi$ as
\be{psidef}
\psi = \int_{0}^{T} ds\, \big(\cS_s \phi\big)
\ee
with $T$ a large constant (suppressed from the notation), implying that
\be{Gap1.71}
-\cA\psi = \phi- \cS_{T} \phi.
\ee
It was shown in \cite{garholmai06}, Lemma 4.3.1, that
\be{martdef}
N_t =\exp\left[\frac{1}{\kappa}\big[\psi(Z_t)-\psi(Z_0)\big]
-\int_0^t ds\, \left(e^{-\frac{1}{\kappa}\psi}\cA e^{\frac{1}{\kappa}\psi}\right)
(Z_s)\right]
\ee
is an exponential $\P_{\eta,x}$-martingale for all $(\eta,x)\in\Omega\times\Z^3$.
Moreover, if we define $\P_{\eta,x}^\new$ in such a way that
\be{Pnewdef}
\P_{\eta,x}^\new(A) = \E_{\,\eta,x}\big(N_t\, \one_A\big)
\ee
for all events $A$ in the $\sigma$-algebra generated by $(Z_s)_{s \in [0,t]}$, then
under $\P_{\eta,x}^\new$ indeed $(Z_s)_{s\geq 0}$ is a Markov process with generator
$\cA^\new$. Using (\ref{Gap1.71}--\ref{Pnewdef}) and
$\E_{\,\nu_\rho,0}^\new=\int_\Omega\nu_\rho(d\eta)\,\E_{\,\eta,0}^\new$, it then follows
that the expectation in (\ref{lscal}) can be written in the form
\be{Gap1.47}
\begin{aligned}
&\E_{\,\nu_\rho,0}\left(\exp\left[\frac{1}{\kappa}\int_0^tds\,
\phi(Z_s)\right]\right)\\
&\qquad= \E_{\,\nu_\rho,0}^\new\bigg(\exp\bigg[
\frac{1}{\kappa}\big[\psi(Z_0) - \psi(Z_t)\big]
+\int_0^t ds\, \left[
\left(e^{-\frac{1}{\kappa}\psi} \cA e^{\frac{1}{\kappa}\psi}\right)
-\cA\left(\frac{1}{\kappa}\psi\right)\right](Z_s)\\
&\qquad\qquad\qquad+\frac{1}{\kappa}\int_0^t ds\, \big(\cS_{T}\phi\big)(Z_s)
\bigg]\bigg).
\end{aligned}
\ee
The first term in the exponent in the r.h.s.\ of (\ref{Gap1.47}) stays bounded as $t\to\infty$
and can therefore be discarded when computing $\lambda^*(\kappa)$ via (\ref{lscal}). We will
see later that the second term and the third term lead to the Green term and the polaron term
in (\ref{limscalfin}), respectively. These terms may be separated from each other with the
help of H\"older's inequality, as stated in Proposition~\ref{lbholder} below.


\subsection{Key propositions}
\label{S2.2}

\bp{lbholder}
For any $\kappa>0$,
\be{Gap1.27}
\lambda^*(\kappa)
\,\,{\leq \atop \geq}\,\,I_1^{q}(\kappa) + I_2^{r}(\kappa)
\ee
with
\be{I12defs}
\begin{aligned}
I_1^{q}(\kappa) &=
\frac1{q}\lim_{t\to\infty}\frac{1}{t}\log\E_{\,\nu_\rho,0}^\new\bigg(\exp\bigg[
q \int_0^t ds\, \left[
\left(e^{-\frac{1}{\kappa}\psi} \cA e^{\frac{1}{\kappa}\psi}\right)
-\cA\left(\frac{1}{\kappa}\psi\right)\right](Z_s) \bigg]\bigg),\\
I_2^{r}(\kappa) &=
\frac1{r}\lim_{t\to\infty}\frac{1}{t}\log\E_{\,\nu_\rho,0}^\new\bigg(\exp\bigg[
\frac{r}{\kappa}\int_0^t ds\, \big(\cS_{T}\phi\big)(Z_s)\bigg]\bigg),
\end{aligned}
\ee
where $1/q + 1/r=1$, with $q>0$, $r>1$ in the first inequality and $q<0$, $0<r<1$
in the second inequality.
\ep

\bpr
See \cite{garholmai06}, Proposition 4.4.1. The existence and finiteness of the limits
in (\ref{I12defs}) follow from Lemma~\ref{varnewlem-1} below.
\epr

By choosing $r$ arbitrarily close to $1$, we see that the proof of our main statement in
(\ref{limscalfin}) reduces to the following two propositions, where we abbreviate
\be{limorders}
\limsup_{t,\kappa,T\to\infty}
=\limsup_{T\to\infty}\limsup_{\kappa\to\infty}\limsup_{t\to\infty}
\quad\text{and}\quad
\lim_{t,\kappa,T\to\infty}
=\lim_{T\to\infty}\lim_{\kappa\to\infty}\lim_{t\to\infty}.
\ee
In the next proposition we write $\psi_T$ instead of $\psi$ to indicate the
dependence on the parameter~$T$.
\bp{mainLem2}
For any $\alpha\in\R$,
\be{main1-3}
\limsup_{t,\kappa,T\to\infty}\frac{\kappa^2}{t}
\log\E_{\,\nu_\rho,0}^\new\left(\exp\left[\alpha\,\int_0^t ds\,
\left[\left(e^{-\frac{1}{\kappa} \psi_T} \cA e^{\frac{1}{\kappa}\psi_T}\right)
-\cA\Big(\frac{1}{\kappa}\psi_T\Big)\right](Z_s)\right]\right)
\leq\frac{\alpha}{6}\,\rho(1-\rho)G.
\ee
\ep

\bp{mainLem1}
For any $\alpha>0$,
\be{remain1-3}
\lim_{t,\kappa,T\to\infty}
\frac{\kappa^2}{t}\log \E_{\,\nu_\rho,0}^\new\bigg(\exp\bigg[
\frac{\alpha}{\kappa}\int_0^t ds\,
\big(\cS_{T}\phi\big)(Z_s)\bigg]\bigg)
= [6\alpha^2\rho(1-\rho)]^2\cP_3.
\ee
\ep

\noindent
These propositions will be proved in Sections~\ref{S3}--\ref{S5}.


\subsection{Preparatory lemmas}
\label{S2.3}

This section contains three elementary lemmas that will be used frequently in
Sections~\ref{S3}--\ref{S5}.

Let $p_t^{(1)}(x,y)$ and $p_t(x,y)=p_t^{(3)}(x,y)$ be the transition kernels of simple
random walk in $d=1$ and $d=3$, respectively, with step rate $1$.

\bl{approxLaplem}
There exists $C>0$ such that, for all $t\geq 0$ and $x,y,e\in\Z^3$ with $\|e\|=1$,
\be{approxLap1}
p_t^{(1)}(x,y) \leq \frac{C}{(1+t)^{\frac12}},
\quad
p_t(x,y)\leq \frac{C}{(1+t)^{\frac{3}{2}}},
\quad
\big|p_t(x+e,y)-p_t(x,y)\big|\leq\frac{C}{(1+t)^2}.
\ee
\el

\bpr
Standard.
\epr

\noindent
(In the sequel we will frequently write $p_t(x-y)$ instead of $p_t(x,y)$.)

From the graphical representation for SSE (Liggett~\cite{lig85}, Chapter VIII, Theorem 1.1)
it is immediate that
\be{graph}
\E_{\,\eta}\big(\xi_t(x)\big) = \sum_{y\in\Z^d}  p_t(x,y)\,\eta(y).
\ee
Recalling (\ref{phidef}--\ref{Ztdef}) and (\ref{psidef}), we therefore have
\be{Gap1.65}
\begin{aligned}
\cS_s \phi (\eta,x)
=\E_{\,\eta,x} \Big(\phi(Z_s)\Big)
&=\E_{\,\eta}\Bigg(\sum_{y\in\Z^3}p_{6s}(x,y)\big[\xi_{s/\kappa}(y)-\rho\big]\Bigg)\\
&=\sum_{z\in\Z^3}p_{6s\onek}(x,z)\big[\eta(z)-\rho\big]
\end{aligned}
\ee
and
\be{Gap1.69}
\psi(\eta,x) = \int_0^{T} ds\, \sum_{z\in\Z^3}
p_{6s\onek}(x,z)\big[\eta(z)-\rho\big],
\ee
where we abbreviate
\be{onekdef}
1[\kappa] = 1 + \frac{1}{6\kappa}.
\ee

\bl{psiboundlem}
For all $\kappa,T>0$, $\eta\in\Omega$, $a,b\in\Z^3$ with $\|a-b\|=1$ and $x\in\Z^3$,
\be{psibound1}
|\psi(\eta,b)-\psi(\eta,a)|\leq 2C\sqrt{T}
\quad\text{for } T\geq1,
\ee
\be{psibound2}
\Big|\psi\big(\eta^{a,b},x\big)-\psi(\eta,x)\Big|\leq 2G,
\ee
\be{psibound3}
\sum_{\{a,b\}}\Big(\psi\big(\eta^{a,b},x\big)-\psi(\eta,x)\Big)^2
\leq \frac{1}{6}G,
\ee
where $C>0$ is the same constant as in Lemma~{\rm \ref{approxLaplem}}, and $G$ is the
value at $0$ of the Green function of simple random walk on $\Z^3$.
\el

\bpr
For a proof of (\ref{psibound2}--\ref{psibound3}), see \cite{garholmai06}, Lemma 4.5.1.
To prove (\ref{psibound1}), we may without loss of generality consider $b=a+e_1$ with
$e_1=(1,0,0)$. Then, by (\ref{Gap1.69}), we have
\be{psibound5}
\begin{aligned}
|\psi(\eta,b)-\psi(\eta,a)|
&\leq \int_0^{T}ds\, \sum_{z\in\Z^3}\big|p_{6s\onek}(z+e_1)-p_{6s\onek}(z)\big|\\
&= \int_0^{T}ds\, \sum_{z\in\Z^3}
\Big|p_{6s\onek}^{(1)}(z_1+e_1)-p_{6s\onek}^{(1)}(z_1)\Big|\, p_{6s\onek}^{(1)}(z_2)\,
p_{6s\onek}^{(1)}(z_3)\\
&= \int_0^{T}ds\, \sum_{z_1\in\Z}
\Big|p_{6s\onek}^{(1)}(z_1+e_1)-p_{6s\onek}^{(1)}(z_1)\Big|\\
&= 2\int_0^{T}ds\, p_{6s\onek}^{(1)}(0)\leq 2C\sqrt{T}.
\end{aligned}
\ee
In the last line we have used the first inequality in (\ref{approxLap1}).
\epr

Let $\cG$ be the Green operator acting on functions $V\colon \Z^3\ra [0,\infty)$ as
\be{greenop}
\cG V(x) = \sum_{y\in\Z^3} G(x-y)V(y),\quad x\in\Z^3,
\ee
with $G(z)=\int_0^\infty dt\, p_t(z)$. Let $\|\cdot\|_\infty$ denote the supremum norm.

\bl{greenoplem}
For all $V\colon \Z^3\ra [0,\infty)$ and $x\in\Z^3$,
\be{greenopeq}
\ES_x^X\bigg(\exp\bigg[\int_0^\infty dt\, V(X_t)\bigg]\bigg)
\leq \Big(1-\|\cG V\|_\infty\Big)^{-1}
\leq \exp\Bigg(\frac{\|\cG V\|_\infty}{1-\|\cG V\|_\infty}\Bigg),
\ee
provided that
\be{greenopasum}
\|\cG V\|_\infty <1.
\ee
\el

\bpr
See \cite{garhol04}, Lemma 8.1.
\epr


\section{Reduction to the original measure}
\label{S3}

In this section we show that the expectations in Propositions~\ref{mainLem2}--\ref{mainLem1}
w.r.t.\ the new measure $\P_{\nu_\rho,0}^\new$ are asymptotically the same as the
expectations w.r.t.\ the old measure $\P_{\nu_\rho,0}$. In Section~\ref{S3.1} we state
a Rayleigh-Ritz formula from which we draw the desired comparison. In Section~\ref{S3.2}
we state the analogues of Propositions~\ref{mainLem2}--\ref{mainLem1} whose proof will
be the subject of Sections~\ref{S4}--\ref{S5}.


\subsection{Rayleigh-Ritz formula}
\label{S3.1}

Recall the definition of $\psi$ in (\ref{psidef}). Let $m$ denote the counting measure
on $\Z^3$. It is easily checked that both $\mu_\rho=\nu_\rho\otimes m$ and $\mu_\rho^\new$
given by
\be{new-dens}
d\mu_\rho^\new = e^{\frac2\kappa\psi}\, d\mu_\rho
\ee
are reversible invariant measures of the Markov processes with generators $\cA$ defined in
(\ref{lb2.7}), respectively, $\cA^\new$ defined in (\ref{lb2i.2}). In particular, $\cA$
and $\cA^\new$ are self-adjoint operators in $L^2(\mu_\rho)$ and $L^2(\mu_\rho^\new)$. Let
$\cD(\cA)$ and $\cD(\cA^\new)$ denote their domains.

\bl{varnewlem-1}
For all bounded measurable $V\colon\Omega\times\Z^3\ra\R$,
\be{varnew-5}
\lim_{t\ra\infty}\frac1t \log
\E_{\,\nu_\rho,0}^\new\bigg(\exp\bigg[\int_0^t ds\,V(Z_s)\bigg]\bigg)
= \sup_{{F\in\cD(\cA^\new)}\atop{\|F\|_{L^2(\mu_\rho^\new)}=1}}
\iint_{\Omega\times\Z^3} d\mu_\rho^\new\,
\Big(V F^2+ F\,\cA^\new F \Big).
\ee
The same is true when $\E_{\,\nu_\rho,0}^\new$, $\mu_\rho^\new$, $\cA^\new$ are
replaced by $\E_{\,\nu_\rho,0}$, $\mu_\rho$, $\cA$, respectively.
\el

\bpr
The limit in the l.h.s.\ of (\ref{varnew-5}) coincides with the upper boundary of the
spectrum of the operator $\cA^\new+V$ on $L^2(\mu_\rho^\new)$, which may be represented
by the Rayleigh-Ritz formula. The latter coincides with the expression in the r.h.s.\ of
(\ref{varnew-5}). The details are similar to \cite{garholmai06}, Section 2.2.
\epr

Lemma~\ref{varnewlem-1} can be used to express the limits as $t\to\infty$ in
Propositions~\ref{mainLem2}--\ref{mainLem1} as variational expressions involving the
new measure. Lemma~\ref{varnewlem-2} below says that, for large $\kappa$, these
variational expressions are close to the corresponding variational expressions for the
old measure. Using Lemma~\ref{varnewlem-1} for the original measure, we may therefore
arrive at the corresponding limit for the old measure.

For later use, in the statement of Lemma~\ref{varnewlem-2} we do not assume that $\psi$
is given by (\ref{psidef}). Instead, we only suppose that $\eta\mapsto \psi(\eta)$ is
bounded and measurable and that there is a constant $K>0$ such
that for all $\eta\in\Omega$, $a,b\in\Z^3$ with $\|a-b\|=1$ and $x\in\Z^3$,
\be{psicond}
|\psi(\eta,b)-\psi(\eta,a)|\leq K
\quad\text{and}\quad
\Big|\psi\big(\eta^{a,b},x\big)-\psi(\eta,x)\Big|\leq K,
\ee
but retain that $\cA^\new$ and $\mu_\rho^\new$ are given by (\ref{lb2i.2}) and (\ref{new-dens}),
respectively.

\bl{varnewlem-2}
Assume {\rm(\ref{psicond})}. Then, for all bounded measurable $V\colon\Omega\times\Z^3\ra\R$,
\be{varnewold-3}
\begin{aligned}
&\sup_{{F\in\cD(\cA^\new)}\atop{\|F\|_{L^2(\mu_\rho^\new)}=1}}
\iint_{\Omega\times\Z^3} d\mu_\rho^\new\,
\Big(V F^2+ F\, \cA^\new F \Big)\\
&\qquad \,\,{\leq \atop \geq}\,\,
e^{\mp\frac{K}{\kappa}}\sup_{{F\in\cD(\cA)}\atop{\|F\|_{L^2(\mu_\rho)}=1}}
\iint_{\Omega\times\Z^3} d\mu_\rho\,
\Big( e^{\pm \frac{K}\kappa} V F^2 + F\, \cA F \Big),
\end{aligned}
\ee
where $\pm$ means $+$ in the first inequality and $-$ in the second inequality,
and $\mp$ means the reverse.
\el

\bpr
Combining (\ref{dL}), (\ref{expro1}) and (\ref{lb2.7}--\ref{lb2i.2}), we have for all
$(\eta,x)\in\Omega\times\Z^3$ and all $F\in\cD(\cA^\new)$,
\be{varnewold-5}
\begin{aligned}
\Big(V F^2+ F\, \cA^\new F \Big)(\eta,x)
&=V(\eta,x)\, F^2(\eta,x)\\
&\quad+ \frac1{6\kappa}\sum_{\{a,b\}}F(\eta,x)\,
e^{\frac1\kappa[\psi(\eta^{a,b},x)-\psi(\eta,x)]}
\Big[F(\eta^{a,b},x)-F(\eta,x)\Big]\\
&\quad+ \sum_{y\colon \|y-x\|=1}F(\eta,x)\,
e^{\frac1\kappa[\psi(\eta,y)-\psi(\eta,x)]}
\big[F(\eta,y)-F(\eta,x)\big].
\end{aligned}
\ee
Therefore, taking into account (\ref{lb2i.2}), (\ref{new-dens}) and the exchangeability
of $\nu_\rho$, we find that
\be{varnewold-7}
\begin{aligned}
\iint_{\Omega\times\Z^3} d\mu_\rho^\new\,
\Big(V F^2+ F\, \cA^\new F \Big)
&=\iint_{\Omega\times\Z^3} d\mu_\rho^\new(\eta,x)\,
\Bigg(V(\eta,x)\, F^2(\eta,x)\\
&\quad- \frac1{12\kappa}\sum_{\{a,b\}}e^{\frac1\kappa[\psi(\eta^{a,b},x)
-\psi(\eta,x)]}\Big[F(\eta^{a,b},x)-F(\eta,x)\Big]^2\\
&\quad- \frac1{2}\sum_{y\colon \|y-x\|=1}e^{\frac1\kappa[\psi(\eta,y)-\psi(\eta,x)]}
\big[F(\eta,y)-F(\eta,x)\big]^2
\Bigg).
\end{aligned}
\ee
Let $\widetilde F=e^{\psi/\kappa}F$. Then, by (\ref{new-dens}) and (\ref{psicond}),
\be{varnewold-9}
\begin{aligned}
(\ref{varnewold-7})
&\,\,{\leq \atop \geq}\,\,
\iint_{\Omega\times\Z^3} d\mu_\rho^\new(\eta,x)\,
\Bigg(V(\eta,x)\, F^2(\eta,x)\\
&\qquad- \frac{e^{\mp \frac{K}\kappa}}{12\kappa}\, \sum_{\{a,b\}}
\Big[F(\eta^{a,b},x)-F(\eta,x)\Big]^2
- \frac{e^{\mp \frac{K}\kappa}}{2}\, \sum_{y\colon \|y-x\|=1}
\big[F(\eta,y)-F(\eta,x)\big]^2 \Bigg)\\
&\,\,=\,\,
\iint_{\Omega\times\Z^3} d\mu_\rho(\eta,x)\,
\Bigg(V(\eta,x)\, \widetilde{F}^{\,2}(\eta,x)\\
&\qquad- \frac{e^{\mp \frac{K}\kappa}}{12\kappa}\, \sum_{\{a,b\}}
\Big[\widetilde F(\eta^{a,b},x)-\widetilde F(\eta,x)\Big]^2
- \frac{e^{\mp \frac{K}\kappa}}{2}\, \sum_{y\colon \|y-x\|=1}
\Big[\widetilde F(\eta,y)-\widetilde F(\eta,x)\Big]^2 \Bigg)\\
&\,\,=\,\,
e^{\mp \frac{K}\kappa}\iint_{\Omega\times\Z^3} d\mu_\rho\,
\Big( e^{\pm \frac{K}\kappa}\, V \widetilde F^2 +\widetilde F\, \cA \widetilde F
\Big).
\end{aligned}
\ee
Taking further into account that
\be{varnewold-8}
\left\|\widetilde F\,\right\|_{L^2(\mu_\rho)}^2
=\left\|F\right\|_{L^2(\mu_\rho^\new)}^2,
\ee
and that $\widetilde F\in\cD(\cA)$ if and only if $F\in\cD(\cA^\new)$, we get the claim.
\epr


\subsection{Reduced key propositions}
\label{S3.2}

At this point we may combine the assertions in Lemmas \ref{varnewlem-1}--\ref{varnewlem-2}
for the potentials
\be{pot1}
V=\alpha\left[\left(e^{-\frac{1}{\kappa} \psi} \cA e^{\frac{1}{\kappa}\psi}\right)
-\cA\Big(\frac{1}{\kappa}\psi\Big)\right]
\ee
and
\be{pot2}
V=\frac\alpha\kappa\big(\cS_{T}\phi\big)
\ee
with $\psi$ given by (\ref{psidef}). Because of (\ref{psibound1}--\ref{psibound2}),
the constant $K$ in (\ref{psicond}) may be chosen to be the maximum of $2G$ and
$2C\sqrt{T}$, resulting in $K/\kappa\to 0$ as $\kappa\to\infty$. Moreover, from
(\ref{psibound3}) and a Taylor expansion of the r.h.s.\ of (\ref{pot1}) we see
that the potential in (\ref{pot1}) is bounded for each $\kappa$ and $T$, and the
same is obviously true for the potential in (\ref{pot2}) because of (\ref{phidef}).
In this way, using a moment inequality to replace the factor $e^{\pm K/\kappa}\alpha$
by a slightly larger, respectively, smaller factor $\alpha^\prime$ independent of $T$
and $\kappa$, we see that the limits in Propositions~\ref{mainLem2}--\ref{mainLem1}
do not change when we replace $\E_{\,\nu_\rho,0}^\new$ by $\E_{\,\nu_\rho,0}$.
Hence it will be enough to prove the following two propositions.

\bp{mainLem2*}
For all $\alpha\in\R$,
\be{main1-3*}
\limsup_{t,\kappa,T\to\infty}\frac{\kappa^2}{t}
\log\E_{\,\nu_\rho,0}\left(\exp\left[\alpha\,\int_0^t ds\,
\left[\left(e^{-\frac{1}{\kappa} \psi} \cA e^{\frac{1}{\kappa}\psi}\right)
-\cA\Big(\frac{1}{\kappa}\psi\Big)\right](Z_s)\right]\right)
\leq \frac{\alpha}{6}\,\rho(1-\rho)G.
\ee
\ep

\bp{mainLem1*}
For all $\alpha>0$,
\be{remain1-3*}
\lim_{t,\kappa,T\to\infty}
\frac{\kappa^2}{t}\log \E_{\,\nu_\rho,0}\bigg(\exp\bigg[
\frac{\alpha}{\kappa}\int_0^t ds\,
\big(\cS_{T}\phi\big)(Z_s)\bigg]\bigg)
= \big[6\alpha^2\rho(1-\rho)\big]^2\cP_3.
\ee
\ep

Proposition \ref{mainLem2*} has already been proven in \cite{garholmai06}, Proposition 4.4.2.
Sections \ref{S4}--\ref{S5} are dedicated to the proof of the lower, respectively, upper
bound in Proposition \ref{mainLem1*}.


\section{Proof of Proposition \ref{mainLem1*}: lower bound}
\label{S4}

In this section we derive the lower bound in Proposition \ref{mainLem1*}.
We fix $\alpha,\kappa,T>0$ and use Lemma \ref{varnewlem-1}, to obtain
\be{lwst1-3}
\begin{aligned}
\lim_{t\to\infty}\frac1t\log\E_{\,\nu_\rho,0}\bigg(\exp\bigg[
\frac\alpha\kappa\int_0^t ds\, \big(\cS_{T}\phi\big)(Z_s)\bigg]\bigg)
=\sup_{{F\in\cD(\cA)}\atop{\|F\|_{L^2(\mu_\rho)}=1}}\iint_{\Omega\times\Z^3}d\mu_\rho
\Big(\frac{\alpha}{\kappa} \big(\cS_{T}\phi\big)F^2
+F \cA F\Big).
\end{aligned}
\ee
In Section~\ref{S4.1} we choose a test function. In Section~\ref{S4.2} we compute and
estimate the resulting expression. In Section~\ref{S4.3} we take the limit $\kappa,T
\to\infty$ and show that this gives the desired lower bound.


\subsection{Choice of test function}
\label{S4.1}

To get the desired lower bound, we use test functions $F$ of the form
\be{testF}
F(\eta,x)=F_1(\eta)F_2(x).
\ee
Before specifying $F_1$ and $F_2$, we introduce some further notation. In addition to
the counting measure $m$ on $\Z^3$, consider the discrete Lebesgue measure $m_\kappa$
on $\Z_\kappa^3=\kappa^{-1}\Z^3$ giving weight $\kappa^{-3}$ to each site in $\Z_\kappa^3$.
Let $l^2(\Z^3)$ and $l^2(\Z_\kappa^3)$ denote the corresponding $l^2$-spaces. Let
$\Delta_\kappa$ denote the lattice Laplacian on $\Z_\kappa^3$ defined by
\be{scal11}
\big(\Delta_\kappa f\big)(x)
=\kappa^2\sum_{{y\in\Z_\kappa^3} \atop {\|y-x\|=\kappa^{-1}}}
\big[f(y)-f(x)\big].
\ee
Choose $f\in\cC_\c^\infty(\R^3)$ with $\|f\|_{L^2(\R^3)}=1$ arbitrarily, where
$\cC_\c^\infty(\R^3)$ is the set of infinitely differentiable functions on $\R^3$
with compact support. Define
\be{scal7}
f_\kappa(x)=\kappa^{-3/2}f\big(\kappa^{-1}x\big),
\quad x\in\Z^3,
\ee
and note that
\be{scal9}
\|f_\kappa\|_{l^2(\Z^3)}=\|f\|_{l^2(\Z_\kappa^3)}
\to 1
\quad\text{as }\kappa\to\infty.
\ee

For $F_2$ choose
\be{F2def}
F_2=\|f_\kappa\|_{l^2(\Z^3)}^{-1}\, f_\kappa.
\ee
To choose $F_1$, introduce the function
\be{phitildef}
\widetilde\phi(\eta)=\frac{\alpha}{\|f_\kappa\|_{l^2(\Z^3)}^2}
\sum_{x\in\Z^3}\big(\cS_{T}\phi\big)(\eta,x)\, f_\kappa^2(x).
\ee
Given $K>0$, abbreviate
\be{STdef}
S=6T\onek
\quad\text{and}\quad
U=6K\kappa^2\onek
\ee
(recall (\ref{onekdef})). For $\kappa>\sqrt{T/K}$, define $\widetilde\psi\colon\Omega\ra
\R$ by
\be{psitildef}
\widetilde\psi=\int_0^{U-S} ds\, \cT_s\widetilde\phi,
\ee
where $(\cT_t)_{t\geq 0}$ is the semigroup generated by the operator $L$ in (\ref{expro1}).
Note that the construction of $\widetilde\psi$ from $\widetilde\phi$ in (\ref{psitildef})
is similar to the construction of $\psi$ from $\phi$ in (\ref{psidef}). In particular,
\be{tildephipsieq}
-L\widetilde\psi = \widetilde\phi-
\cT_{U-S} \widetilde\phi.
\ee
Combining the probabilistic representations of the semigroups $(\cS_t)_{t\geq 0}$
(generated by $\cA$ in (\ref{lb2.7})) and $(\cT_t)_{t\geq 0}$ (generated by $L$ in
(\ref{expro1})) with the graphical representation formulas (\ref{graph}--\ref{Gap1.65}),
and using (\ref{scal7}--\ref{scal9}), we find that
\be{phitildef*}
\widetilde\phi(\eta)
=\frac{\alpha}{\|f\|_{l^2(\Z^3_\kappa)}^2}
\int_{\Z_\kappa^3}m_\kappa(dx)\,
f^2(x)\sum_{z\in\Z^3}p_{S}(\kappa x,z)[\eta(z)-\rho]
\ee
and
\be{psitildef*}
\widetilde\psi(\eta)
=\sum_{z\in\Z^3} h(z)[\eta(z)-\rho]
\ee
with
\be{psitildef**}
h(z)=\frac{\alpha}{\|f\|_{l^2(\Z^3_\kappa)}^2}
\int_{\Z_\kappa^3} m_\kappa(dx)\, f^2(x)\int_{S}^{U}ds\,p_{s}(\kappa x,z).
\ee
Using the second inequality in (\ref{approxLap1}), we have
\be{hbd}
0\leq h(z)\leq \frac{C\alpha}{\sqrt{T}},\quad z\in\Z^3.
\ee
Now choose $F_1$ as
\be{F1def}
F_1=\big\|e^{\widetilde\psi}\big\|_{L^2(\nu_\rho)}^{-1}\,
e^{\widetilde\psi}.
\ee

For the above choice of $F_1$ and $F_2$, we have $\|F_1\|_{L^2(\nu_\rho)}=\|F_2\|_{l^2(\Z^3)}=1$
and, consequently, $\|F\|_{L^2(\mu_\rho)}=1$. With $F_1$, $F_2$ and $\widetilde\phi$ as above,
and $\cA$ as in (\ref{lb2.7}), after scaling space by $\kappa$ we arrive at the following lemma.

\bl{testfunct}
For $F$ as in {\rm (\ref{testF})}, {\rm (\ref{F2def})} and {\rm (\ref{F1def})},
all $\alpha,T,K>0$ and $\kappa>\sqrt{T/K}$,
\be{lwst2-13}
\begin{aligned}
&\kappa^2\iint_{\Omega\times\Z^3} d\mu_\rho\,\Big(
\frac\alpha\kappa\big(\cS_{T}\phi\big)F^2+F\cA F\Big)\\
&\qquad= \frac{1}{\|f\|_{l^2(\Z_\kappa^3)}^2}\int_{\Z_\kappa^3}d m_\kappa\, f \Delta_\kappa f
+\frac{\kappa}{\|e^{\widetilde\psi}\|_{L^2(\nu_\rho)}^2}
\int_\Omega d\nu_\rho\Big(\widetilde\phi e^{2\widetilde\psi}
+e^{\widetilde\psi} Le^{\widetilde\psi}\Big),
\end{aligned}
\ee
where $\widetilde \phi$ and $\widetilde \psi$ are as in {\rm (\ref{phitildef})} and
{\rm (\ref{psitildef})}.
\el


\subsection{Computation of the r.h.s.\ of (\ref{lwst2-13})}
\label{S4.2}

Clearly, as $\kappa\to\infty$ the first summand in the r.h.s.\ of (\ref{lwst2-13})
converges to
\be{lwst3-3}
\int_{\R^3} dx\, f(x)\, \Delta f(x)
=-\big\|\nabla_{\R^3}f\big\|_{L^2(\R^3)}^2.
\ee
The computation of the second summand in the r.h.s.\ of (\ref{lwst2-13}) is more
delicate:

\bl{2summlem}
For all $\alpha>0$ and $0<\epsilon<K$,
\be{2summbd}
\begin{aligned}
&\liminf_{\kappa,T\to\infty}
\frac{\kappa}{\|e^{\widetilde\psi}\|_{L^2(\nu_\rho)}^2}
\int_\Omega d\nu_\rho\Big(\widetilde\phi e^{2\widetilde\psi}
+e^{\widetilde\psi} Le^{\widetilde\psi}\Big)\\
&\qquad\geq 6\alpha^2\rho(1-\rho)\int_{\R^3} dx\, f^2(x) \int_{\R^3} dy\, f^2(y)
\Bigg(\int_{6\epsilon}^{6K} dt\, p_t^{(G)}(x,y)
-\int_{6K}^{12K} dt\, p_t^{(G)}(x,y)\Bigg),
\end{aligned}
\ee
where
\be{Gaussianker}
p_t^{(G)}(x,y)=(4\pi t)^{-3/2}\exp[-\|x-y\|^2/4t]
\ee
denotes the Gaussian transition kernel associated with $\Delta_{\R^3}$, the continuous
Laplacian on $\R^3$.
\el

\bpr
Using the probability measure
\be{lwst4-3}
d\nu_\rho^\new= \big\|e^{\widetilde\psi}\big\|_{L^2(\nu_\rho)}^{-2}\,
e^{2\widetilde\psi}\, d\nu_\rho
\ee
in combination with (\ref{tildephipsieq}), we may write the term under the $\liminf$ in
(\ref{2summbd}) in the form
\be{lwst4-3*}
\kappa\int_\Omega d\nu_\rho^\new\, \Big(e^{-\widetilde\psi}
L e^{\widetilde\psi}-L\widetilde\psi +\cT_{U-S}\widetilde\phi\Big).
\ee
This expression can be handled by making a Taylor expansion of the $L$-terms and showing
that the $\cT_{U-S}$-term is nonnegative. Indeed, by the definition of $L$ in (\ref{expro1}),
we have
\be{lwst5-3}
\Big(e^{-\widetilde\psi} L e^{\widetilde\psi}-L\widetilde\psi\Big)(\eta)
=\frac1{6}\sum_{\{a,b\}}\bigg(e^{[\widetilde\psi(\eta^{a,b})
-\widetilde\psi(\eta)]}
-1-\Big[\widetilde\psi\big(\eta^{a,b}\big)-\widetilde\psi(\eta)\Big]\bigg).
\ee
Recalling the expressions for $\widetilde\psi$ in (\ref{psitildef*}--\ref{psitildef**}) and
using (\ref{hbd}), we get for $a,b\in\Z^3$ with $\|a-b\|=1$,
\be{lwst5-5}
\big|\widetilde\psi\big(\eta^{a,b}\big)-\widetilde\psi(\eta)\big|
=|h(a)-h(b)|\, |\eta(b)-\eta(a)|\\
\leq\frac{C\alpha}{\sqrt{T}}.
\ee
Hence, a Taylor expansion of the exponent in the r.h.s.\ of (\ref{lwst5-3}) gives
\be{lwst5-7}
\int_\Omega d\nu_\rho^\new
\Big(e^{-\widetilde\psi}\, L e^{\widetilde\psi}-L\widetilde\psi\Big)
\geq\frac{e^{-C\alpha/\sqrt{T}}}{12}\int_\Omega d\nu_\rho^\new
\sum_{\{a,b\}}\Big[\widetilde\psi\big(\eta^{a,b}\big)
-\widetilde\psi(\eta)\Big]^2.
\ee
Using (\ref{psitildef*}), we obtain
\be{lwst5-11}
\int_\Omega\nu_\rho^\new (d\eta)\sum_{\{a,b\}}
\Big[\widetilde\psi\big(\eta^{a,b}\big)-\widetilde\psi(\eta)\Big]^2
=\sum_{\{a,b\}}\big[h(a)-h(b)\big]^2\int_\Omega\nu_\rho^\new (d\eta)
\big[\eta(b)-\eta(a)\big]^2.
\ee
Using (\ref{lwst4-3}), we have (after cancellation of factors not depending on $a$ or $b$)
\be{lwst5-13}
\int_\Omega\nu_\rho^\new (d\eta)\big[\eta(b)-\eta(a)\big]^2
=\frac{\displaystyle \int_\Omega\nu_\rho(d\eta)\,
e^{2\chi_{a,b}(\eta)}\big[\eta(b)-\eta(a)\big]^2}
{\displaystyle \int_\Omega\nu_\rho(d\eta)\, e^{2\chi_{a,b}(\eta)}}
\ee
with
\be{lwst5-15}
\chi_{a,b}(\eta)=h(a)\eta(a)+h(b)\eta(b).
\ee
Using (\ref{hbd}), we obtain that
\be{lwst5-18}
\int_\Omega\nu_\rho^\new (d\eta)\big[\eta(b)-\eta(a)\big]^2
\geq e^{-4C\alpha/\sqrt{T}}\int_\Omega\nu_\rho(d\eta)\big[\eta(b)-\eta(a)\big]^2
=e^{-4C\alpha/\sqrt{T}}2\rho(1-\rho).
\ee
On the other hand, by (\ref{psitildef**}),
\be{lwst5-19}
\begin{aligned}
\sum_{\{a,b\}}\big[h(a)-h(b)\big]^2
&=\frac{\alpha^2}{\|f\|_{l^2(\Z^3_\kappa)}^4}
\int_{S}^{U}dt\int_{S}^{U}ds\int_{\Z_\kappa^3} m_\kappa(dx)\,
f^2(x)\int_{\Z_\kappa^3} m_\kappa(dy)\, f^2(y)\\
&\qquad\qquad\times\sum_{\{a,b\}}\big[p_t(\kappa x,a)-p_t(\kappa x,b)\big]
\big[p_s(\kappa y,a)-p_s(\kappa y,b)\big]
\end{aligned}
\ee
with
\be{lwst5-21}
\begin{aligned}
\sum_{\{a,b\}}\big[p_t(\kappa x,a)-p_t(\kappa x,b)\big]
\big[p_s(\kappa y,a)-p_s(\kappa y,b)\big]
&= -\sum_{a\in\Z^3}p_t(\kappa x,a)
\Delta p_s(\kappa x,a)\\
&=-6\sum_{a\in\Z^3}p_t(\kappa x,a)
\bigg(\frac{\partial}{\partial s}p_s(\kappa y,a)\bigg),
\end{aligned}
\ee
where $\Delta$ acts on the first spatial variable of $p_s(\cdot\,,\cdot)$ and
$\Delta p_s=6(\partial p_s/\partial s)$. Therefore,
\be{lwst5-23}
\begin{aligned}
(\ref{lwst5-19})
&=6\int_{S}^{U}dt \int_{\Z_\kappa^3} m_\kappa(dx)\,
f^2(x)\int_{\Z_\kappa^3} m_\kappa(dy)\, f^2(y)\sum_{a\in\Z^3}p_t(\kappa x,a)
\big[p_{S}(\kappa y,a)-p_{U}(\kappa y,a)\big]\\
&=6\int_{\Z_\kappa^3} m_\kappa(dx)\, f^2(x)\int_{\Z_\kappa^3} m_\kappa(dy)\, f^2(y)
\Bigg(\int_{2S}^{S+U} dt\, p_t(\kappa x,\kappa y)
-\int_{U+S}^{2U} dt\, p_{t}(\kappa x,\kappa y)\Bigg).
\end{aligned}
\ee
Combining (\ref{lwst5-7}--\ref{lwst5-11}) and (\ref{lwst5-18}--\ref{lwst5-19})
and (\ref{lwst5-23}), we arrive at
\be{lwst5-24}
\begin{aligned}
\int_\Omega d\nu_\rho^\new
\Big(e^{-\widetilde\psi}L e^{\widetilde\psi}-L\widetilde \psi\Big)
&\geq \frac{e^{-5C\alpha/\sqrt{T}}\alpha^2}{\|f\|_{l^2(\Z^3_\kappa)}^4}
\rho(1-\rho)\int_{\Z_\kappa^3}
m_\kappa(dx)\, f^2(x) \int_{\Z_\kappa^3} m_\kappa(dy)\, f^2(y)\\
&\qquad\qquad\times\Bigg(\int_{2S}^{S+U} dt\, p_t(\kappa x,\kappa y)
-\int_{U+S}^{2U} dt\, p_{t}(\kappa x,\kappa y)\Bigg).
\end{aligned}
\ee
After replacing $2S$ in the first integral by $6\epsilon\kappa^2\onek$, using a
Gaussian approximation of the transition kernel $p_t(x,y)$ and recalling the
definitions of $S$ and $U$ in (\ref{STdef}), we get that, for any $\epsilon>0$,
\be{lwst5-25}
\begin{aligned}
&\liminf_{\kappa,T\to\infty} \kappa\int_\Omega d\nu_\rho^\new
\Big(e^{-\widetilde\psi}L e^{\widetilde\psi}-L\widetilde \psi\Big)\\
&\qquad\geq 6\alpha^2\rho(1-\rho)\int_{\R^3} dx\, f^2(x) \int_{\R^3} dy\, f^2(y)
\Bigg(\int_{6\epsilon}^{6K} dt\, p_t^{(G)}(x,y)
-\int_{6K}^{12K} dt\, p_t^{(G)}(x,y)\Bigg).
\end{aligned}
\ee

At this point it only remains to check that the $\cT_{U-S}$-term in (\ref{lwst4-3*}) is
nonnegative. By (\ref{phitildef*}) and the probabilistic representation of the semigroup
$(\cT_t)_{t\geq 0}$, we have
\be{lwst4-11}
\int_\Omega d\nu_\rho^\new\, \cT_{U-S}\widetilde\phi
=\frac{\alpha}{\|f\|_{l^2(\Z_\kappa^3)}^2}
\int_{\Z_\kappa^3} m_\kappa(dx)\, f^2(x) \sum_{z\in\Z^3} p_{U}(\kappa x,z)
\int_\Omega \nu_\rho^\new(d\eta)[\eta(z)-\rho]
\ee
and, by (\ref{lwst4-3}),
\be{lwst4-15}
\begin{aligned}
\int_\Omega \nu_\rho^\new(d\eta) [\eta(z)-\rho]
&=-\rho+\frac{\rho e^{2h(z)}}{\rho e^{2h(z)}+1-\rho}
=-\rho+\frac{\rho}{1-(1-\rho)\big(1- e^{-2h(z)}\big)}\\
&\geq -\rho+\rho\Big[1+(1-\rho)\Big(1-e^{-2h(z)}\Big)\Big]
=\rho(1-\rho)\Big(1-e^{-2h(z)}\Big),
\end{aligned}
\ee
which proves the claim.
\epr


\subsection{Proof of the lower bound in Proposition~\ref{mainLem1*}}
\label{S4.3}

We finish by using Lemma \ref{2summlem} to prove the lower bound in
Proposition~\ref{mainLem1*}.

\bpr
Combining (\ref{lwst2-13}--\ref{2summbd}), we get
\be{lwst4-16}
\begin{aligned}
&\liminf_{\kappa,T\to\infty}
\kappa^2\iint_{\Omega\times\Z^3} d\mu_\rho\,\Big(
\frac\alpha\kappa\big(\cS_{T}\phi\big)F^2-F\cA F\Big)\\
&\qquad\geq 6\alpha^2\rho(1-\rho)\int_{\R^3} dx\, f^2(x) \int_{\R^3} dy\, f^2(y)
\Bigg(\int_{6\epsilon}^{6K} dt\, p_t^{(G)}(x,y)
-\int_{6K}^{12K} dt\, p_t^{(G)}(x,y)\Bigg)\\
&\qquad\quad-\big\|\nabla_{\R^3}f\big\|_{L^2(\R^3)}^2.
\end{aligned}
\ee
Letting $\epsilon\da 0$, $K\to\infty$, replacing $f(x)$ by $\gamma^{3/2}f(\gamma x)$ with
$\gamma=6\alpha^2\rho(1-\rho)$, taking the supremum over all $f\in C_\c^\infty(\R^3)$
such that $\|f\|_{L^2(\R^3)}=1$ and recalling (\ref{lwst1-3}), we arrive at
\be{lwbdineq}
\liminf_{t,\kappa,T\to\infty}
\frac{\kappa^2}{t}\log \E_{\,\nu_\rho,0}\bigg(\exp\bigg[
\frac{\alpha}{\kappa}\int_0^t ds\,
\big(\cS_{T}\phi\big)(Z_s)\bigg]\bigg)
\geq \big[6\alpha^2\rho(1-\rho)\big]^2\cP_3,
\ee
which is the desired inequality.
\epr


\section{Proof of Proposition \ref{mainLem1*}: upper bound}
\label{S5}

In this section we prove the upper bound in Proposition~\ref{mainLem1*}. The proof is
long and technical. In Sections~\ref{S5.1} we ``freeze'' and ``defreeze''
the exclusion dynamics on long time intervals. This allows us to approximate the relevant
functionals of the random walk in terms of its occupation time measures on those intervals.
In Section~\ref{S5.defreez*} we use a spectral bound to reduce the study of the long-time
asymptotics for the resulting time-dependent potentials to the investigation of
time-independent potentials. In Section~\ref{S5.3} we make a cut-off for small times,
showing that these times are negligible in the limit as $\kappa\to\infty$, perform a
space-time scaling and compactification of the underlying random walk, and apply a large
deviation principle for the occupation time measures, culminating in the appearance of
the variational expression for the polaron term $\cP_3$.


\subsection{Freezing, defreezing and reduction to two key lemmas}
\label{S5.1}


\subsubsection{Freezing}
\label{S5.1.1}

We begin by deriving a preliminary upper bound for the expectation
in Proposition \ref{mainLem1*} given by
\be{upst1-3}
\E_{\,\nu_\rho,0}\bigg(\exp\bigg[\int_0^t ds\, V(Z_s)\bigg]\bigg)
\ee
with
\be{upst2-3}
V(\eta,x)
=\frac\alpha\kappa\big(\cS_{T}\phi\big)(\eta,x)
=\frac\alpha\kappa\sum_{y\in\Z^3}p_{6T\onek}(x,y)(\eta(y)-\rho),
\ee
where, as before, $T$ is a large constant. To this end, we divide the time interval
$[0,t]$ into $\lfloor t/R_\kappa\rfloor$ intervals of length
\be{Rkappadef}
R_\kappa=R\kappa^2
\ee
with $R$ a large constant, and ``freeze'' the exclusion dynamics $(\xi_{t/\kappa})_{t\geq 0}$
on each of these intervals. As will become clear later on, this procedure allows us to
express the dependence of (\ref{upst1-3}) on the random walk $X$ in terms of objects that
are close to integrals over occupation time measures of $X$ on time intervals of length
$R_\kappa$. We will see that the resulting expression can be estimated from above by
``defreezing'' the exclusion dynamics. We will subsequently see that, after we have taken
the limits $t\to\infty$, $\kappa\to\infty$ and $T\to\infty$, the resulting estimate can
be handled by applying a large deviation principle for the space-time rescaled occupation
time measures in the limit as $R\to\infty$. The latter will lead us to the polaron term.

Ignoring the negligible final time interval $[\lfloor t/R_\kappa\rfloor R_\kappa,t]$,
using H\"older's inequality with $p,q>1$ and $1/p+1/q=1$, and inserting (\ref{upst2-3}),
we see that (\ref{upst1-3}) may be estimated from above as
\be{upst1-11}
\begin{aligned}
&\E_{\,\nu_\rho,0}\bigg(
\exp\bigg[\int_0^{\lfloor t/R_\kappa \rfloor R_\kappa} ds\,V(Z_s)\bigg]\bigg)\\
&\quad=\E_{\,\nu_\rho,0}\Bigg(\exp\Bigg[\frac\alpha\kappa\sum_{k=1}^{\lfloor t/R_\kappa \rfloor}
\int_{(k-1)R_\kappa}^{k R_\kappa} ds \sum_{y\in\Z^3} p_{6T\onek}(X_s,y)
\big(\xi_{s/\kappa}(y)-\rho\big)\Bigg]\Bigg)\\
&\quad\leq \Big(\cE_{R,\alpha q}^{(1)}(t)\Big)^{1/q}
\Big(\cE_{R,\alpha p}^{(2)}(t)\Big)^{1/p}
\end{aligned}
\ee
with
\be{E1def}
\begin{aligned}
\cE_{R,\alpha}^{(1)}(t)=\cE_{R,\alpha}^{(1)}(\kappa,T;t)
&=\E_{\,\nu_\rho,0}\Bigg(\exp\Bigg[\frac{\alpha}\kappa
\sum_{k=1}^{\lfloor t/R_\kappa \rfloor}
\int_{(k-1)R_\kappa}^{k R_\kappa} ds \sum_{y\in\Z^3}\bigg(p_{6T\onek}(X_s,y)\,
\xi_{\frac{s}\kappa}(y)\\
&\qquad\qquad\qquad\quad
-p_{6T\onek+\frac{s-(k-1)R_\kappa}{\kappa}}(X_s,y)\, \xi_{\frac{(k-1)R_\kappa}\kappa}(y)
\bigg)\Bigg]\Bigg)
\end{aligned}
\ee
and
\be{E2def}
\begin{aligned}
&\cE_{R,\alpha}^{(2)}(t)=\cE_{R,\alpha}^{(2)}(\kappa,T;t)\\
&\quad= \E_{\,\nu_\rho,0}\Bigg(\exp\Bigg[\frac{\alpha}\kappa
\sum_{k=1}^{\lfloor t/R_\kappa \rfloor}
\int_{(k-1)R_\kappa}^{k R_\kappa} ds \sum_{y\in\Z^3}
p_{6T\onek+\frac{s-(k-1)R_\kappa}{\kappa}}(X_s,y)
\Big(\xi_{\frac{(k-1)R_\kappa}\kappa}(y)-\rho\Big)\Bigg]\Bigg).
\end{aligned}
\ee
Therefore, by choosing $p$ close to $1$, the proof of the upper bound in
Proposition \ref{mainLem1*} reduces to the proof of the following two lemmas.

\bl{freez}
For all $R,\alpha>0$,
\be{freez1}
\limsup_{t,\kappa,T\to\infty}\frac{\kappa^2}{t}\log \cE_{R,\alpha}^{(1)}(\kappa,T;t)\leq 0.
\ee
\el

\bl{freez*}
For all $\alpha>0$,
\be{freez*1}
\limsup_{R\to\infty}\limsup_{t,\kappa,T\to\infty}\frac{\kappa^2}{t}\log
\cE_{R,\alpha}^{(2)}(\kappa,T;t)
\leq \big[6\alpha^2\rho(1-\rho)\big]^2\cP_3.
\ee
\el

Lemma~\ref{freez} will be proved in Section~\ref{S5.freez}, Lemma~\ref{freez*} in
Sections~\ref{S5.freez*}--\ref{S5.4}.


\subsubsection{Proof of Lemma~\ref{freez}}
\label{S5.freez}

\bpr
Fix $R,\alpha>0$ arbitrarily. Given a path $X$, an initial configuration $\eta\in\Omega$
and $k\in\N$, we first derive an upper bound for
\be{frzeq3}
\E_{\,\eta}\Bigg(\exp\Bigg[\frac{\alpha}\kappa\int_{0}^{R_\kappa} ds \sum_{y\in\Z^3}
\bigg(p_{6T\onek}\Big(X_s^{(k,\kappa)},y\Big)\xi_{\frac{s}\kappa}(y)
-p_{6T\onek+\frac{s}{\kappa}}\Big(X_s^{(k,\kappa)},y\Big) \eta(y)\bigg)\Bigg]\Bigg),
\ee
where
\be{frzeq5}
X_s^{(k,\kappa)}=X_{(k-1)R_\kappa+s}.
\ee
To this end, we use the independent random walk approximation $\widetilde\xi$ of $\xi$
(cf.\ \cite{garholmai06}, Proposition 1.2.1), to obtain
\be{frzeq7}
(\ref{frzeq3})
\leq \prod_{y\in A_\eta}\ES_{\,0}^{Y}\Bigg(\exp\Bigg[\frac{\alpha}\kappa\int_{0}^{R_\kappa} ds
\bigg(p_{6T\onek}\Big(X_s^{(k,\kappa)},y+Y_{\frac{s}\kappa}\Big)
-p_{6T\onek+\frac{s}{\kappa}}\Big(X_s^{(k,\kappa)},y\Big)\bigg)\Bigg]\Bigg),
\ee
where $Y$ is simple random walk on $\Z^3$ with jump rate $1$ (i.e., with generator
$\frac16\Delta$), $\ES_0^Y$ is expectation w.r.t.\ $Y$ starting from $0$, and
\be{frzeq9}
A_\eta=\{x\in\Z^3\colon \eta(x)=1\}.
\ee
Observe that the expectation w.r.t.\ $Y$ of the expression in the exponent is zero. Therefore,
a Taylor expansion of the exponential function yields the bound
\be{frzeq11}
\begin{aligned}
&\ES_{\,0}^{Y}\Bigg(\exp\Bigg[\frac{\alpha}\kappa\int_{0}^{R_\kappa} ds
\bigg(p_{6T\onek}\Big(X_s^{(k,\kappa)},y+Y_{\frac{s}\kappa}\Big)
-p_{6T\onek+\frac{s}{\kappa}}\Big(X_s^{(k,\kappa)},y\Big)\bigg)\Bigg]\Bigg)\\
&\quad\leq 1+\sum_{n=2}^{\infty}\prod_{l=1}^{n}\Bigg(\frac{\alpha}{\kappa}
\int_{s_{l-1}}^{R_{\kappa}} ds_{l} \sum_{y_l\in\Z^3}
p_{\frac{s_l-s_{l-1}}{\kappa}}(y_{l-1},y_l)\\
&\qquad\qquad\qquad\times\bigg[
p_{6T\onek}\Big(X_{s_l}^{(k,\kappa)},y+y_l\Big)
+p_{6T\onek+\frac{s_l}{\kappa}}\Big(X_{s_l}^{(k,\kappa)},y\Big)\bigg]\Bigg),
\end{aligned}
\ee
where $s_0=0$, $y_0=0$, and the product has to be understood in a noncommutative way.
Using the Chapman-Kolmogorov equation and the inequality $p_t(z)\leq p_t(0)$,
$z\in\Z^3$, we find that
\be{frzeq13}
\begin{aligned}
&\int_{s_{l-1}}^{R_{\kappa}} ds_{l} \sum_{y_l\in\Z^3}
p_{\frac{s_l-s_{l-1}}{\kappa}}(y_{l-1},y_l)\,\,
\bigg[
p_{6T\onek}\Big(X_{s_l}^{(k,\kappa)},y+y_l\Big)
+p_{6T\onek+\frac{s_l}{\kappa}}\Big(X_{s_l}^{(k,\kappa)},y\Big)\bigg]\Bigg)\\
&\qquad \leq 2\int_0^\infty ds\, p_{T+\frac{s}{\kappa}}(0)=2\kappa G_T(0)
\end{aligned}
\ee
with
\be{frzeq15}
G_T(0)=\int_T^\infty ds\, p_s(0)
\ee
the cut-off Green function of simple random walk at $0$ at time $T$. Substituting this into
the above bound for $l=n,n-1,\cdots,3$, computing the resulting geometric series, and using
the inequality $1+x\leq e^x$, we obtain
\be{frzeq17}
\begin{aligned}
(\ref{frzeq11})
&\leq \exp\Bigg[\frac{C_T\alpha^2}{\kappa^2}\prod_{l=1}^{2}\int_{s_{l-1}}^{R_\kappa} ds_l
\sum_{y_l\in\Z^3} p_{\frac{s_l-s_{l-1}}{\kappa}}(y_{l-1},y_l)\\
&\qquad\qquad\times\bigg(
p_{6T\onek}\Big(X_{s_l}^{(k,\kappa)},y+y_l\Big)
+p_{6T\onek+\frac{s_l}{\kappa}}\Big(X_{s_l}^{(k,\kappa)},y\Big)\bigg)\Bigg]
\end{aligned}
\ee
with
\be{frzeq19}
C_T=\frac{1}{1-2\alpha G_T(0)},
\ee
provided that $2\alpha G_T(0)<1$, which is true for $T$ large enough. Note that
$C_T\to 1$ as $T\to\infty$. Substituting (\ref{frzeq17}) into (\ref{frzeq7}), we
find that
\be{frzeq21}
\begin{aligned}
(\ref{frzeq3})
&\leq\exp\Bigg[\frac{C_T\alpha^2}{\kappa^2}\sum_{y\in\Z^3}\prod_{l=1}^{2}
\int_{s_{l-1}}^{R_\kappa} ds_l \sum_{y_l\in\Z^3}
p_{\frac{s_l-s_{l-1}}{\kappa}}(y_{l-1},y_l)\\
&\qquad\qquad\times\bigg(
p_{6T\onek}\Big(X_{s_l}^{(k,\kappa)},y+y_l\Big)
+p_{6T\onek+\frac{s_l}{\kappa}}\Big(X_{s_l}^{(k,\kappa)},y\Big)\bigg)\Bigg].
\end{aligned}
\ee
Using once more the Chapman-Kolmogorov equation and $p_t(x,y)=p_t(x-y)$, we may compute
the sums in the exponent, to arrive at
\be{frzeq23}
\begin{aligned}
(\ref{frzeq3})
\leq\exp\Bigg[&\frac{C_T\alpha^2}{\kappa^2}\int_{0}^{R_\kappa}
ds_1 \int_{s_1}^{R_\kappa} ds_2
\bigg(p_{12T\onek+\frac{s_2-s_1}{\kappa}}
\Big(X_{s_2}^{(k,\kappa)}-X_{s_1}^{(k,\kappa)}\Big)\\
&\quad
+3p_{12T\onek+\frac{s_2+s_1}{\kappa}}\Big(X_{s_2}^{(k,\kappa)}
-X_{s_1}^{(k,\kappa)}\Big)\bigg)\Bigg].
\end{aligned}
\ee
Note that this bound does not depend on the initial configuration $\eta$ and depends on
the process $X$ only via its increments on the time interval $[(k-1)R_\kappa,kR_\kappa]$.
By (\ref{frzeq5}), the increments over the time intervals labelled $k=1,2,\cdots,\lfloor
t/R_\kappa\rfloor$ are independent and identically distributed. Using $\E_{\nu_\rho,0}
=\int \nu_\rho(d\eta)\ES_0^X\E_\eta$, we can therefore apply the Markov property of the
exclusion dynamics $(\xi_{t/\kappa})_{t\geq 0}$ at times $R_\kappa$, $2R_\kappa$, $\cdots$,
($\lfloor t/R_\kappa\rfloor-1)R_\kappa$ to the expectation in the r.h.s.\ of (\ref{E1def}),
insert the bound (\ref{frzeq23}) and afterwards use that $(X_t)_{t\geq 0}$ has independent
increments, to arrive at
\be{frzeq25}
\begin{aligned}
\log \cE_{R,\alpha}^{(1)}(t)
\leq \frac{t}{R_\kappa}\log \ES_{\,0}^X\Bigg(\exp\Bigg[&\frac{C_T\alpha^2}{\kappa^2}
\int_{0}^{R_\kappa} ds_1 \int_{s_1}^{R_\kappa} ds_2
\bigg(p_{12T\onek+\frac{s_2-s_1}{\kappa}}\big(X_{s_2}-X_{s_1}\big)\\
&\quad+3p_{12T\onek+\frac{s_2+s_1}{\kappa}}\big(X_{s_2}-X_{s_1}\big)\bigg)\Bigg]\Bigg).
\end{aligned}
\ee
Hence, recalling the definition of $R_\kappa$ in (\ref{Rkappadef}), we obtain
\be{frzeq27}
\begin{aligned}
&\limsup_{t\to\infty} \frac{\kappa^2}{t}\log \cE_{R,\alpha}^{(1)}(t)\\
&\quad\leq\frac1R\log \ES_{\,0}^X
\Bigg(\exp\Bigg[\frac{C_T\alpha^2R}{R_\kappa}\int_{0}^{R_\kappa} ds_1
\int_{s_1}^{R_\kappa} ds_2
\bigg(p_{12T\onek+\frac{s_2-s_1}{\kappa}}\big(X_{s_2}-X_{s_1}\big)\\
&\qquad\qquad\qquad\qquad\qquad
+3p_{12T\onek+\frac{s_2+s_1}{\kappa}}\big(X_{s_2}-X_{s_1}\big)\bigg)\Bigg]\Bigg).
\end{aligned}
\ee

Let
\be{Xhat}
\widehat X_t=X_t+Y_{t/\kappa},
\ee
and let $\ES_0^{\widehat X}=\ES_0^X \ES_0^Y$ be the expectation
w.r.t.\ $\widehat X$ starting at $0$. Observe that
\be{frzeq29}
p_{t+s/\kappa}(z) = \ES_0^Y\Big(p_t\big(z+Y_{s/\kappa}\big)\Big).
\ee
We next apply Jensen's inequality w.r.t.\ the first integral in the r.h.s.\ of (\ref{frzeq27}),
substitute $s_2=s_1+s$, take into account that $X$ has independent increments, and afterwards
apply Jensen's inequality w.r.t.\ $\ES_0^Y$, to arrive at the following upper bound for
the expectation in (\ref{frzeq27}):
\be{}
\begin{aligned}
&\ES_{\,0}^X\Bigg(\exp\Bigg[\frac{C_T\alpha^2R}{R_\kappa}\int_{0}^{R_\kappa} ds_1
\int_{s_1}^{R_\kappa} ds_2
\bigg(p_{12T\onek+\frac{s_2-s_1}{\kappa}}\big(X_{s_2}-X_{s_1}\big)\\
&\qquad\qquad\quad
+3p_{12T\onek+\frac{s_2+s_1}{\kappa}}\big(X_{s_2}-X_{s_1}\big)\bigg)\Bigg]\\
&\quad\leq \frac{1}{R_\kappa}\int_0^{R_\kappa} ds_1\, \ES_0^X\Bigg(\exp\Bigg[
C_T\alpha^2R\int_0^\infty ds\, \ES_0^Y\bigg(p_{12T\onek}\Big(X_{s}+Y_{\frac{s}{\kappa}}\Big)\\
&\qquad\qquad\qquad\qquad\qquad\qquad\quad
+3p_{12T\onek+\frac{2s_1}{\kappa}}\Big(X_{s}+Y_{\frac{s}{\kappa}}\Big)\bigg)\Bigg]\Bigg)\\
&\quad\leq \frac{1}{R_\kappa}\int_0^{R_\kappa} ds_1\, \ES_0^{\widehat X}\Bigg(\exp\Bigg[
C_T\alpha^2R\int_0^\infty ds\Big(p_{12T\onek}\big(\widehat X_{s}\big)
+3p_{12T\onek+\frac{2s_1}{\kappa}}\big(\widehat X_{s}\big)\Big)\Bigg]\Bigg).
\end{aligned}
\ee
Applying Lemma~\ref{greenoplem}, we can bound the last expression from above by
\be{Ghatbd}
\exp\Bigg[\frac{4C_T\alpha^2R \widehat G_{2T}(0)}{1-4C_T\alpha^2R \widehat G_{2T}(0)}\Bigg],
\ee
where $\widehat G_{2T}(0)$ is the cut-off at time $2T$ of the Green function $\widehat G$ at
$0$ for $\widehat X$ (which has generator $\onek \Delta$). Since $\widehat G_{2T}(0)\to
\frac16 G_{12T}(0)$ as $\kappa\to\infty$, and since the latter converges to zero as
$T\to\infty$, a combination of the above estimates with (\ref{frzeq27}) gives the claim.
\epr


\subsubsection{Defreezing}
\label{S5.freez*}

To prove Lemma \ref{freez*}, we next ``defreeze'' the exclusion dynamics in $\cE_{R,
\alpha}^{(2)}(t)$. This can be done in a similar way as the ``freezing'' we did in
Section~\ref{S5.1.1}, by taking into account the following remarks. In (\ref{E2def}),
each single summand is asymptotically negligible as $t\to\infty$. Hence, we can safely
remove a summand at the beginning and add a summand at the end. After that we can
bound the resulting expression from above with the help of H\"older's inequality with
weights $p,q>1$, $1/p+1/q=1$, namely,
\be{Holderagain}
\begin{aligned}
&\E_{\,\nu_\rho,0}\Bigg(\exp\Bigg[\frac{\alpha}\kappa
\sum_{k=1}^{\lfloor t/R_\kappa \rfloor}
\int_{k R_\kappa}^{(k+1) R_\kappa} ds \sum_{y\in\Z^3}
p_{6T\onek+\frac{s-k R_\kappa}{\kappa}}(X_s,y)
\Big(\xi_{\frac{k R_\kappa}\kappa}(y)-\rho\Big)\Bigg]\Bigg)\\
&\quad\leq\Big(\cE_{R,\alpha q}^{(3)}(t)\Big)^{1/q}
\Big(\cE_{R,\alpha p}^{(4)}(t)\Big)^{1/p}
\end{aligned}
\ee
with
\be{E3def}
\begin{aligned}
&\cE_{R,\alpha}^{(3)}(t)=\cE_{R,\alpha}^{(3)}(\kappa,T;t)\\
&\quad=\E_{\,\nu_\rho,0}\Bigg(\exp\Bigg[\frac{\alpha}{\kappa R_\kappa}
\sum_{k=1}^{\lfloor t/R_\kappa \rfloor}
\int_{(k-1)R_\kappa}^{k R_\kappa} du \int_{kR_\kappa}^{(k+1)R_\kappa} ds
\sum_{y\in\Z^3}\bigg(p_{6T\onek+\frac{s-kR_\kappa}{\kappa}}(X_s,y)\,
\xi_{\frac{kR_\kappa}\kappa}(y)\\
&\qquad\qquad\qquad\qquad
-p_{6T\onek+\frac{s-u}{\kappa}}(X_s,y)\, \xi_{\frac{u}\kappa}(y)
\bigg)\Bigg]\Bigg)
\end{aligned}
\ee
and
\be{E4def}
\begin{aligned}
&\cE_{R,\alpha}^{(4)}(t)=\cE_{R,\alpha}^{(4)}(\kappa,T;t)\\
&\;= \E_{\,\nu_\rho,0}\Bigg(\exp\Bigg[\frac{\alpha}{\kappa R_\kappa}
\sum_{k=1}^{\lfloor t/R_\kappa \rfloor}
\int_{(k-1)R_\kappa}^{k R_\kappa} du \int_{k R_\kappa}^{(k+1)R_\kappa}ds
\sum_{y\in\Z^3} p_{6T\onek+\frac{s-u}{\kappa}}(X_s,y)
\Big(\xi_{\frac{u}\kappa}(y)-\rho\Big)\Bigg]\Bigg).
\end{aligned}
\ee
In this way, choosing $p$ close to $1$, we see that the proof of Lemma \ref{freez*}
reduces to the proof of the following two lemmas.

\bl{defreez}
For all $R,\alpha>0$,
\be{E3-1}
\limsup_{t,\kappa,T\to\infty}\frac{\kappa^2}{t}\log\cE_{R,\alpha}^{(3)}(\kappa,T;t)\leq 0.
\ee
\el
\bl{defreez*}
For all $\alpha>0$,
\be{E4-1}
\limsup_{R\to\infty}\limsup_{t,\kappa,T\to\infty}\frac{\kappa^2}{t}
\log\cE_{R,\alpha}^{(4)}(\kappa,T;t)
\leq \big[6\alpha^2\rho(1-\rho)\big]^2\cP_3.
\ee
\el

In the remaining sections we prove Lemmas~\ref{defreez}--\ref{defreez*} and thereby complete
the proof of the upper bound in Proposition~\ref{mainLem1*}.


\subsubsection{Proof of Lemma~\ref{defreez}}
\label{S5.defreez}

\bpr
The proof goes along the same lines as the proof of Lemma \ref{freez}.
Instead of (\ref{frzeq3}), we consider
\be{defreest1}
\begin{aligned}
\E_{\,\eta}\Bigg(\exp\Bigg[
&\frac{\alpha}{\kappa R_\kappa}\int_{0}^{R_\kappa} du
\int_{R_\kappa}^{2R_\kappa} ds \sum_{y\in\Z^3}
\bigg(p_{6T\onek+\frac{s-R_\kappa}{\kappa}}\Big(X_s^{(k,\kappa)},y\Big)
\xi_{\frac{R_\kappa}\kappa}(y)\\
&-p_{6T\onek+\frac{s-u}{\kappa}}\Big(X_s^{(k,\kappa)},y\Big)
\xi_{\frac{u}{\kappa}}(y)\bigg)\Bigg]\Bigg).
\end{aligned}
\ee
Applying Jensen's inequality w.r.t.\ the first integral and the Markov property
of the exclusion dynamics $(\xi_{t/\kappa})_{t\geq 0}$ at time $u/\kappa$, we see
that it is enough to derive an appropriate upper bound for
\be{Ezetabd}
\begin{aligned}
\E_{\,\zeta}\Bigg(\exp\Bigg[
&\frac{\alpha}{\kappa}\int_{R_\kappa}^{2R_\kappa} ds \sum_{y\in\Z^3}
\bigg(p_{6T\onek+\frac{s-R_\kappa}{\kappa}}\Big(X_s^{(k,\kappa)},y\Big)
\xi_{\frac{R_\kappa-u}\kappa}(y)\\
&-p_{6T\onek+\frac{s-u}{\kappa}}\Big(X_s^{(k,\kappa)},y\Big) \zeta(y)\bigg)\Bigg]\Bigg)
\end{aligned}
\ee
uniformly in $\zeta\in\Omega$ and $u\in[0,R_\kappa]$. The main steps are the
same as in the proof of Lemma \ref{freez}. Instead of (\ref{frzeq23}), we obtain
\be{}
\begin{aligned}
(\ref{Ezetabd})
\leq\exp\Bigg[&\frac{C_T\alpha^2}{\kappa^2}
\int_{R_\kappa}^{2R_\kappa} ds_1 \int_{s_1}^{2R_\kappa} ds_2
\bigg(p_{12T\onek+\frac{s_2-s_1}{\kappa}+\frac{2(s_1-R_\kappa)}{\kappa}}
\Big(X_{s_2}^{(k,\kappa)}-X_{s_1}^{(k,\kappa)}\Big)\\
&\quad
+3p_{12T\onek+\frac{s_2-s_1}{\kappa}+\frac{2(s_1-u)}{\kappa}}
\Big(X_{s_2}^{(k,\kappa)}-X_{s_1}^{(k,\kappa)}\Big)\bigg)\Bigg],
\end{aligned}
\ee
and this expression may be bounded from above by (\ref{Ghatbd}).
\epr

\subsection{Spectral bound}
\label{S5.defreez*}

The advantage of Lemma~\ref{defreez*} compared to the original upper bound in
Proposition~\ref{mainLem1*} is that, modulo a small time correction of the form
$(s-u)/\kappa$, the expression under the expectation in (\ref{E4def}) depends on
$X$ only via its occupation time measures on the time intervals $[kR_\kappa,(k+1)R_\kappa]$,
$k=1,2,\cdots,\lfloor t/R_\kappa\rfloor$. This will allow us in Section~\ref{S5.3}
to use a large deviation principle for these occupation time measures. The present
section consists of five steps, organized in Sections \ref{S5.2}--\ref{S5.reddefreez*},
leading up to a final lemma that will be proved in Section~\ref{S5.3}.

We abbreviate
\be{Vpot}
V_{k,u}(\eta)
=V_{k,u}^{\kappa,X}(\eta)
=\frac{1}{R_\kappa}\int_{kR_\kappa}^{(k+1)R_\kappa}ds\sum_{y\in\Z^3}
p_{6T\onek+\frac{s-u}{\kappa}}\big(X_s,y\big)(\eta(y)-\rho)
\ee
and rewrite the expression for $\cE_{R,\alpha}^{(4)}(t)$ in (\ref{E4def}) in the form
\be{E4-3}
\cE_{R,\alpha}^{(4)}(t)
=\E_{\,\nu_\rho,0}\Bigg(\exp\Bigg[\frac{\alpha}\kappa\sum_{k=1}^{\lfloor t/R_\kappa \rfloor}
\int_{(k-1)R_\kappa}^{k R_\kappa} du\, V_{k,u}\big(\xi_{u/\kappa}\big)\Bigg]\Bigg).
\ee
In (\ref{Vpot}) and subsequent expressions we suppress the dependence on $T$ and $R$.


\subsubsection{Reduction to a spectral bound}
\label{S5.2}

Let $B(\Omega)$ denote the Banach space of bounded measurable
functions on $\Omega$ equipped with the supremum norm $\|\cdot\|_\infty$. Given
$V\in B(\Omega)$,
let
\be{upst4-3}
\lambda(V)=\lim_{t\ra\infty}\frac1t\log\E_{\,\nu_\rho}\bigg(\exp\bigg[
\int_0^t V(\xi_s)\, ds\bigg]\bigg)
\ee
denote the associated Lyapunov exponent. The limit in (\ref{upst4-3}) exists and coincides
with the upper boundary of the spectrum of the self-adjoint operator $L+V$ on $L^2(\nu_\rho)$,
written
\be{spect7}
\lambda(V)=\sup\Sp(L+V).
\ee

\bl{spectrallem}
For all $t>0$ and all bounded and piecewise continuous $V\colon [0,t]\ra B(\Omega)$,
\be{spect1}
\E_{\,\nu_\rho}\bigg(\exp\bigg[\int_0^t V_u(\xi_u)\, du\bigg]\bigg)
\leq \exp\bigg[\int_0^t \lambda(V_s)\, ds\bigg].
\ee
\el

\bpr
In the proof we will assume that $s\mapsto V_s$ is continuous. The extension to piecewise
continuous $s\mapsto V_s$ will be straightforward. Let $0=t_0<t_1<\cdots<t_r=t$ be a
partition of the interval $[0,t]$. Then
\be{spect5}
\begin{aligned}
\int_0^t V_u(\xi_u)\, du
&\leq\sum_{k=1}^{r}\int_{t_{k-1}}^{t_k} V_{t_{k-1}}(\xi_s)\, ds
+\sum_{k=1}^{r}\max_{s\in[t_{k-1},t_k]}\|V_s-V_{t_{k-1}}\|_\infty
\big(t_k-t_{k-1}\big)\\
&\leq\sum_{k=1}^{r}\int_{t_{k-1}}^{t_k} V_{t_{k-1}}(\xi_s)\, ds
+t\max_{k=1,\cdots,r}\max_{s\in[t_{k-1},t_k]}\|V_s-V_{t_{k-1}}\|_\infty.
\end{aligned}
\ee
Let $(\cS_t^V)_{t\geq 0}$ denote the semigroup generated by $L+V$ on $L^2(\nu_\rho)$
with inner product $(\cdot\,,\cdot)$ and norm $\|\cdot\|$. Then
\be{spect9}
\big\|\cS_t^V\big\|=e^{t\lambda(V)}.
\ee
Using the Markov property, we find that
\be{spect11}
\begin{aligned}
\E_{\,\nu_\rho}\Bigg(\exp\Bigg[\sum_{k=1}^{r}\int_{t_{k-1}}^{t_k}
V_{t_{k-1}}(\xi_s)\, ds\Bigg]\Bigg)
&=\bigg(\cS_{t_1}^{V_{t_0}}\, \cS_{t_2-t_1}^{V_{t_1}}
\cdots \cS_{t_r-t_{r-1}}^{V_{t_{r-1}}}\,\one,\one\bigg)\\
&\leq \big\|\cS_{t_1}^{V_{t_0}}\big\|\, \big\|\cS_{t_2-t_1}^{V_{t_1}}\big\|
\cdots \big\|\cS_{t_r-t_{r-1}}^{V_{t_{r-1}}}\big\|\\
&=\exp\Bigg[\sum_{k=1}^{r}\lambda\big(V_{t_{k-1}}\big)(t_k-t_{k-1})\Bigg].
\end{aligned}
\ee
Combining (\ref{spect5}) and (\ref{spect11}), we arrive at
\be{spect13}
\log\E_{\,\nu_\rho}\bigg(\int_0^t V_s(\xi_s)\, ds\bigg)
\leq \sum_{k=1}^{r}\lambda\big(V_{t_{k-1}}\big)(t_k-t_{k-1})
+t\max_{k=1,\cdots,r}\max_{s\in[t_{k-1},t_k]}\big\|V_s-V_{t_{k-1}}\big\|_\infty.
\ee
Since the map $V\mapsto\lambda(V)$ from $B(\Omega)$ to $\R$ is continuous (which can
be seen e.g.\ from (\ref{spect9}) and the Feynman-Kac representation of $\cS_t^V$),
the claim follows by letting the mesh of the partition tend to zero.
\epr

\bl{spectralbdlem}
For all $\alpha,T,R,t,\kappa>0$,
\be{upst4-11}
\E_{\,\nu_\rho,0}
\Bigg(\exp\Bigg[\frac\alpha\kappa\sum_{k=1}^{\lfloor t/R_\kappa\rfloor}
\int_{(k-1)R_\kappa}^{k R_\kappa} du\,
V_{k,u}\big(\xi_{u/\kappa}\big)\Bigg]\Bigg)
\leq\ES_{\,0}^{X}\Bigg(\exp\Bigg[\sum_{k=1}^{\lfloor t/R_\kappa\rfloor}
\int_{(k-1)R_\kappa}^{k R_\kappa} du\,
\lambda_{k,u} \Bigg]\Bigg)
\ee
with
\be{upst4-13}
\lambda_{k,u}
=\lambda_{k,u}^{\kappa,X}
=\lim_{t\to\infty}\frac1t\log\E_{\,\nu_\rho}
\bigg(\exp\bigg[\frac\alpha\kappa \int_0^t ds\,
V_{k,u}^{\kappa,X}\big(\xi_{s/\kappa}\big)\bigg]\bigg),
\ee
where $u\in[(k-1)R_\kappa,k R_\kappa]$, $k=1,2,\cdots,\lfloor t/R_\kappa\rfloor$.
\el

\bpr
Apply Lemma~\ref{spectrallem} to the potential $V_u(\eta)=(\alpha/\kappa)V_{k,u}(\eta)$
for $u\in[(k-1)R_\kappa,kR_\kappa]$ with $(\xi_u)_{u\geq 0}$ replaced by
$(\xi_{u/\kappa})_{u\geq 0}$, and take the expectation w.r.t.\ $\ES_0^X$.
\epr

The spectral bound in Lemma \ref{spectralbdlem} enables us to estimate the expression
in (\ref{E4-3}) from above by finding upper bounds for the expectation in (\ref{upst4-13})
with a \emph{time-independent} potential $V_{k,u}$. This goes as follows. Fix $\kappa$,
$X$, $k$ and $u$, and abbreviate
\be{upst5-3}
\widehat\phi=\alpha V_{k,u}^{\kappa,X}.
\ee
Let $(\cQ_t)_{t\geq 0}$ be the semigroup generated by $(1/\kappa)L$, and define
\be{psihatdef}
\widehat\psi
=\int_0^M dr\,\big(\cQ_r\widehat\phi\,\big)
\ee
with
\be{Sregime}
M= 3K\onek\kappa^3
\ee
for a large constant $K>0$. Then
\be{psihateq}
-\frac{1}{\kappa} L \widehat\psi = \widehat\phi-\cQ_M \widehat\phi
\ee
with
\be{Qphihat}
\begin{aligned}
\big(\cQ_r \widehat\phi\,\big)(\eta)
&=\frac\alpha{R_\kappa}
\int_{k R_\kappa}^{(k+1)R_\kappa} ds \sum_{y\in\Z^3}
p_{6T\onek+\frac{s-u+r}\kappa}(X_s,y)\big[\eta(y)-\rho\big]\\
&=\alpha\sum_{y\in\Z^3} \Xi_r(y)[\eta(y)-\rho]
\end{aligned}
\ee
and
\be{uplem3-3}
\Xi_r(x) = \Xi_{k,u,r}^{\kappa,X}(x)
=\frac{1}{R_\kappa}\int_{k R_\kappa}^{(k+1)R_\kappa} ds\,
p_{6T\onek+\frac{s-u+r}{\kappa}}(X_s,x).
\ee
As in Section \ref{S2}, we introduce new probability measures $\P_\eta^\new$ by an absolute
continuous transformation of the probability measures $\P_\eta$, in the same way as in
(\ref{martdef}--\ref{Pnewdef}) with $\psi$ and $\cA$ replaced by $\widehat\psi$ and $(1/\kappa) L$,
respectively. Under $\P_\eta^\new$, $(\xi_{t/\kappa})_{t\geq 0}$ is a Markov process with generator
\be{Lnewgen}
\frac1\kappa L^\new f
= e^{-\frac1\kappa \widehat\psi} \frac1\kappa L \bigg(e^{\frac1\kappa \widehat\psi}f\bigg)
-\bigg(e^{-\frac1\kappa \widehat\psi} \frac1\kappa L e^{\frac1\kappa \widehat\psi}\bigg)f.
\ee
Since $\eta\mapsto\widehat\psi(\eta)$ is bounded, we have, similarly as in
Proposition~\ref{lbholder} with $q=r=2$,
\be{upst5-11}
\lambda_{k,u}^{\kappa,X}
\leq \limsup_{t\to\infty}\frac1{2t}\log\Big(\cE^{(5)}_{k,u}(t)\Big)
+\limsup_{t\to\infty}\frac1{2t}\log\Big(\cE^{(6)}_{k,u}(t)\Big)
\ee
with
\be{upst5-13}
\cE^{(5)}_{k,u}(t)=\cE^{(5)}_{k,u}(\kappa,X;t)
=\E_{\,\nu_\rho}^\new\bigg(\exp\bigg[
\frac{2}\kappa \int_0^t dr\, \bigg[
\Big(e^{-\frac{1}{\kappa}\widehat\psi} L e^{\frac{1}{\kappa}\widehat\psi}\,\Big)
-L\bigg(\frac{1}{\kappa}\widehat\psi\bigg)\bigg]\big(\xi_{r/\kappa}\big) \bigg]\bigg)
\ee
and
\be{upst5-13*}
\cE^{(6)}_{k,u}(t)=\cE^{(6)}_{k,u}(\kappa,X;t)
=\E_{\,\nu_\rho}^\new\bigg(\exp\bigg[
\frac{2}{\kappa}\int_0^t dr\, \big(\cQ_{M}\widehat\phi\,\big)\big(\xi_{r/\kappa}\big)
\bigg]\bigg),
\ee
where $\E_{\nu_\rho}^\new=\int_\Omega \nu_\rho(d\eta)\, \E_\eta^\new$, and we suppress the
dependence on the constants $T$, $K$, $R$.


\subsubsection{Two further lemmas}
\label{S5.tfls}

For $a,b\in\Z^3$ with $\|a-b\|=1$, define
\be{uplem3-2}
\begin{aligned}
&\cK_{k,u}(a,b)=\cK_{k,u}^{\kappa,X}(a,b)=
e^{2C\alpha/T}\frac{\alpha^2}{3\kappa^3}\, \int_0^M dr\int_r^M d\widetilde r\,
\big[\Xi_r(a)-\Xi_r(b)\big]\big[\Xi_{\widetilde r}(a)-\Xi_{\widetilde r}(b)\big]
\end{aligned}
\ee
with $\Xi_r$ given by (\ref{uplem3-3}) and $C$ the constant from Lemma \ref{approxLaplem}.
Abbreviate
\be{calKdef}
\big\|\cK_{k,u}\big\|_1=\sum_{\{a,b\}}\cK_{k,u}^{\kappa,X}(a,b).
\ee

\bl{uplem3}
For all $\alpha,T,K,R,\kappa,t>0$, $u\in[(k-1)R_\kappa,k R_\kappa]$,
$k=1,2,\cdots,\lfloor t/R_\kappa\rfloor$, and all paths $X$,
\be{uplem3-1*}
\cE^{(5)}_{k,u}(t)
\leq \E_{\,\nu_\rho}\Bigg(\exp\Bigg[
\kappa \big\|\cK_{k,u}\big\|_1\int_0^{t/\kappa} dr\,
\big[\xi_{r}(e_1)-\xi_{r}(0)\big]^2\Bigg]\Bigg)
\ee
with
\be{K1bd}
\big\|\cK_{k,u}\big\|_1
\leq e^{2C\alpha/T}\frac{2\alpha^2}{\kappa^2 R_\kappa^2}
\int_{k R_\kappa}^{(k+1) R_\kappa}ds \int_{k R_\kappa}^{(k+1) R_\kappa} d\widetilde s
\int_0^M dr\, p_{12T\onek+\frac{s+\widetilde s-2u+2r}\kappa}\big(X_{\widetilde s}-X_s\big).
\ee
\el

\bl{st5lem}
There exists $\kappa_0>0$ such that for all $\kappa>\kappa_0$, $K>1$,
$\alpha,T,R,\kappa,t>0$, $u\in[(k-1)R_\kappa,k R_\kappa]$,
$k=1,2,\cdots,\lfloor t/R_\kappa\rfloor$, and all paths $X$,
\be{st5lem-1}
\cE^{(6)}_{k,u}(t)
\leq \exp\Bigg(\frac{D_{\alpha,T,K}}{\kappa^2}\rho t\Bigg),
\ee
where the constant $D_{\alpha,T,K}$ does not depend on $R$, $t$, $\kappa$, $u$ or $k$
and satisfies
\be{st5lem-2}
\lim_{K\to\infty}D_{\alpha,T,K}=0,
\quad\text{uniformly in } T\geq 1.
\ee
\el

\subsubsection{Proof of Lemma~\ref{uplem3}}
\label{S5.uplem}

\bpr
We want to replace $\E_{\nu_\rho}^\new$ by $\E_{\nu_\rho}$ in formula (\ref{upst5-13})
by applying the analogues of Lemmas \ref{varnewlem-1} and \ref{varnewlem-2}. To this end,
we need to compute the constant $K$ in (\ref{psicond}) for $\psi$ replaced by $\widehat\psi$.
Recalling (\ref{psihatdef}) and (\ref{Qphihat}), we have, for $\eta\in\Omega$ and $a,b\in\Z^3$
with $\|a-b\|=1$,
\be{upst6-7}
\widehat\psi(\eta^{a,b})-\widehat\psi(\eta)
=\alpha \int_0^M dr\,\big[\Xi_r(a)-\Xi_r(b)\big][\eta(b)-\eta(a)].
\ee
Hence,
\be{psihatdif}
\begin{aligned}
\Big|\widehat\psi\big(\eta^{a,b}\big)-\widehat\psi(\eta)\Big|
\leq\alpha \int_0^M dr\,\big|\Xi_r(a)-\Xi_r(b)\big|
\leq C\alpha \int_0^\infty dr\,
\bigg(1+6T+\frac{r}\kappa\bigg)^{-2}
\leq \frac{C\alpha}T \kappa.
\end{aligned}
\ee
Here we have used (\ref{uplem3-3}) and the right-most inequality in (\ref{approxLap1}).
This yields
\be{upst6-2}
\cE^{(5)}_{k,u}(t)
\leq \E_{\,\nu_\rho}\bigg(\exp\bigg[
\frac{2}\kappa e^{C\alpha/T}\int_0^t dr\, \bigg[
\Big(e^{-\frac{1}{\kappa}\widehat\psi} L e^{\frac{1}{\kappa}\widehat\psi}\,\Big)
-L\bigg(\frac{1}{\kappa}\widehat\psi\bigg)\bigg]\big(\xi_{r/\kappa}\big) \bigg]\bigg).
\ee
By (\ref{expro1}), we have
\be{upst6-3}
\frac1\kappa\bigg[e^{-\frac1\kappa\widehat\psi}Le^{\frac1\kappa\widehat\psi}
-L\Big(\frac{1}{\kappa}\widehat\psi\Big)\bigg](\eta)
=\frac{1}{6\kappa}\sum_{\{a,b\}}
\bigg(e^{\frac1\kappa[\widehat\psi(\eta^{a,b})-\widehat\psi(\eta)]}-1
-\frac1\kappa\Big[\widehat\psi(\eta^{a,b})-\widehat\psi(\eta)\Big]\bigg).
\ee
In view of (\ref{psihatdif}), a Taylor expansion of the r.h.s.\ of (\ref{upst6-3}) gives
\be{upst6-5}
\frac1\kappa\bigg[e^{-\frac1\kappa\widehat\psi}Le^{\frac1\kappa\widehat\psi}
-L\Big(\frac{1}{\kappa}\widehat\psi\Big)\bigg](\eta)
\leq\frac{e^{C\alpha/T}}{12\kappa^3}\sum_{\{a,b\}}
\Big(\widehat\psi(\eta^{a,b})-\widehat\psi(\eta)\Big)^2.
\ee
Hence, recalling (\ref{uplem3-2}) and (\ref{upst6-7}), we get
\be{upst6-9}
\begin{aligned}
&\E_{\,\nu_\rho}\bigg(\exp\bigg[
\frac{2}\kappa e^{C\alpha/T}\int_0^t dr\, \bigg[
\Big(e^{-\frac{1}{\kappa}\widehat\psi} L e^{\frac{1}{\kappa}\widehat\psi}\,\Big)
-L\Big(\frac{1}{\kappa}\widehat\psi\Big)\bigg]\big(\xi_{r/\kappa}\big) \bigg]\bigg)\\
&\qquad\leq\E_{\,\nu_\rho}\bigg(\exp\bigg[
\int_0^t dr \sum_{\{a,b\}}\cK_{k,u}(a,b)\Big[\xi_{\frac{r}{\kappa}}(b)
-\xi_{\frac{r}{\kappa}}(a)\Big]^2\bigg]\bigg).
\end{aligned}
\ee
Using Jensen's inequality w.r.t.\ the probability kernel $\cK_{k,u}/\|\cK_{k,u}\|_1$,
together with the translation invariance of $\xi$ under $\P_{\nu_\rho}$, we arrive at
(\ref{uplem3-1*}). To derive (\ref{K1bd}), observe that for arbitrary
$h,\widetilde h, r, \widetilde r>0$ and $x,y\in\Z^d$,
\be{}
\begin{aligned}
&\sum_{\{a,b\}}\Big[p_{h+\frac{r}\kappa}(x,a)-p_{h+\frac{r}\kappa}(x,b)\Big]
\Big[p_{\widetilde h+\frac{\widetilde r}\kappa}(y,a)
-p_{\widetilde h+\frac{\widetilde r}\kappa}(y,b)\Big]\\
&\quad=-\sum_{a\in\Z^3}p_{h+\frac{r}\kappa}(x,a)\Delta p_{\widetilde h
+\frac{\widetilde r}{\kappa}}(y,a)
=-6\kappa\sum_{a\in\Z^3}p_{h+\frac{r}\kappa}(x,a)\,
 \frac{\partial}{\partial \widetilde r}\, p_{\widetilde h+\frac{\widetilde r}{\kappa}}(y,a),
\end{aligned}
\ee
where $\Delta$ acts on the first spatial variable of $p_t(\cdot,\cdot)$ and $\frac16\Delta
p_{t/\kappa}=\kappa(\partial/\partial t)p_{t/\kappa}$. Recalling (\ref{uplem3-3}),
it follows that
\be{}
\sum_{\{a,b\}}\big[\Xi_r(a)-\Xi_r(b)\big]\big[\Xi_{\widetilde r}(a)-\Xi_{\widetilde r}(b)\big]
=-6\kappa \sum_{a\in\Z^3}\Xi_r(a)\frac{\partial}{\partial \widetilde r}\, \Xi_{\widetilde r}(a)
\ee
and, consequently,
\be{upst6-23}
\begin{aligned}
\big\|\cK_{k,u}\big\|_1
&=e^{2C\alpha/T}\frac{2\alpha^2}{\kappa^2}\int_0^M dr \sum_{a\in\Z^3}
\Xi_r(a)\big[\Xi_r(a)-\Xi_M(a)\big]\\
&\leq e^{2C\alpha/T}\frac{2\alpha^2}{\kappa^2}\int_0^M dr \sum_{a\in\Z^3}
\Xi_r(a)^2.
\end{aligned}
\ee
Hence, taking into account (\ref{uplem3-3}), we arrive at (\ref{K1bd}).
\epr

\subsubsection{Proof of Lemma~\ref{st5lem}}
\label{S5.st5lem}

\bpr
Using the same arguments as in (\ref{psihatdif}--\ref{upst6-2}), we can replace
$\E_{\nu_\rho}^\new$ by $\E_{\nu_\rho}$ in formula (\ref{upst5-13*}), to obtain
\be{upst5-21}
\cE^{(6)}_{k,u}(t)
\leq\E_{\,\nu_\rho}\bigg(\exp\bigg[
\frac{2}{\kappa} e^{C\alpha/T}\int_0^t dr\,
\big(\cQ_{M}\widehat\phi\,\big)\big(\xi_{r/\kappa}\big)\bigg]\bigg).
\ee
Because of (\ref{Qphihat}), this yields
\be{upst5-22}
\exp\bigg[\frac{2\alpha}{\kappa}e^{C\alpha/T}\rho t\bigg]\cE^{(6)}_{k,u}(t)\\
\leq\E_{\,\nu_\rho}\bigg(\exp\bigg[
\frac{2\alpha}{\kappa} e^{C\alpha/T}\int_0^t dr\,
\sum_{y\in\Z^3}\Xi_M(y)\xi_{r/\kappa}(y)\bigg]\bigg).
\ee
Now, using the independent random walk approximation $\widetilde \xi$ of $\xi$
(see \cite{garholmai06}, Proposition 1.2.1), we find that
\be{upst5-23}
\begin{aligned}
&\E_{\,\nu_\rho}\bigg(\exp\bigg[
\frac{2\alpha}{\kappa} e^{C\alpha/T}\int_0^t dr\,
\sum_{y\in\Z^3}\Xi_M(y)\xi_{r/\kappa}(y)\bigg]\bigg)\\
&\quad\leq \int\nu_\rho(d\eta)\prod_{x\in A_\eta}\ES_x^Y
\bigg(\exp\bigg[\frac{2\alpha}{\kappa}e^{C\alpha/T}\int_0^t dr\,
\Xi_M\big(Y_{r/\kappa}\big)\bigg]\bigg),
\end{aligned}
\ee
where $A_\eta$ is given by (\ref{frzeq9}) and $Y$ is simple random walk with step rate
$1$. Define
\be{vdef}
v(x,t)=\ES_x^{Y}\bigg(\exp\bigg[
\frac{2\alpha}\kappa e^{C\alpha/T}\int_0^t dr\,
\Xi_M\big(Y_{r/\kappa}\big)\bigg]\bigg),
\quad (x,t)\in\Z^3\times[0,\infty),
\ee
and write
\be{wdef}
w(x,t)=v(x,t)-1.
\ee
Then we may bound (\ref{upst5-22}) from above as follows:
\be{upst5-24}
\begin{aligned}
\text{r.h.s.\ } (\ref{upst5-22})
&\leq\int\nu_\rho(d\eta)\prod_{x\in\Z^3}\big[1+\eta(x)w(x,t)\big]\\
&=\prod_{x\in\Z^3}\big[1+\rho w(x,t)\big]\\
&\leq \exp\bigg(\rho\sum_{x\in\Z^3}w(x,t)\bigg).
\end{aligned}
\ee
By the Feynman-Kac formula, $w$ is the solution of the Cauchy problem
\be{wcauchy}
\frac{\partial}{\partial t}w(x,t)=\frac{1}{6\kappa}\Delta w(x,t)+
\frac{2\alpha}\kappa e^{C\alpha/T}\Xi_M(x)\big[1+w(x,t)\big],
\quad w(\cdot,0)\equiv 0.
\ee
Therefore
\be{wcauchysum}
\frac{\partial}{\partial r}\sum_{x\in\Z^3}w(x,r)
=\frac{2\alpha}\kappa e^{C\alpha/T}\sum_{x\in\Z^3}\Xi_M(x)\big[1+w(x,r)\big].
\ee
Integrating (\ref{wcauchysum}) w.r.t.\ $r$ over the time interval $[0,t]$ and
substituting the resulting expression into (\ref{upst5-24}), we get
\be{wbd}
\begin{aligned}
\text{r.h.s.\ } (\ref{upst5-22})
\leq \exp\Bigg[\frac{2\alpha}\kappa e^{C\alpha/T}\rho
\int_0^t dr\, \sum_{x\in\Z^3}\Xi_M(x) \big(1+w(x,r)\big)\Bigg].
\end{aligned}
\ee
Since $\sum_{x\in\Z^3}\Xi_M(x)=1$, this leads to
\be{upst5-31}
\begin{aligned}
\cE^{(6)}_{k,u}(t)
\leq
\exp\Bigg[\frac{2\alpha}\kappa e^{C\alpha/T}\rho
\int_0^t dr\, \sum_{x\in\Z^3}\Xi_M(x) w(x,r)\Bigg].
\end{aligned}
\ee
An application of Lemma \ref{greenoplem} to the expectation in the r.h.s.\
of (\ref{vdef}) gives
\be{upst5-33}
v(x,t)\leq \Big(1-2\alpha e^{C\alpha/T}\big\|\cG \Xi_M\big\|_\infty\Big)^{-1}.
\ee
Next, using (\ref{Sregime}) and (\ref{uplem3-3}), we find that
\be{upst5-35}
\big\|\cG \Xi_M\big\|_\infty
\leq G_{6T+M/\kappa}(0)\leq G_{3K\kappa^2}(0),
\ee
where the r.h.s.\ tends to zero as $\kappa\to\infty$. Thus, if $K> 1$ and $\kappa>\kappa_0$
with $\kappa_0$ large enough (not depending on the other parameters), then $v(x,t)\leq 2$,
and hence $w(x,t)\leq 1$, for all $x\in\Z^3$ and $t\geq 0$, so that (\ref{wcauchy}) implies
that $w\leq\widehat w$, where $\widehat w$ solves
\be{upst5-47}
\frac{\partial}{\partial t}\widehat w(x,t)=\frac{1}{6\kappa}\Delta \widehat w(x,t)+
\frac{4\alpha}\kappa e^{C\alpha/T} \Xi_M(x),
\qquad \widehat w(\cdot,0)\equiv 0.
\ee
The solution of this Cauchy problem has the representation
\be{upst5-49}
\begin{aligned}
\widehat w(x,t)
=\frac{4\alpha}\kappa e^{C\alpha/T} \int_0^t dr\,\sum_{y\in\Z^3}
p_{\frac{r}{\kappa}}(x,y)\Xi_M(y)
=\frac{4\alpha}\kappa e^{C\alpha/T} \int_0^t dr\, \Xi_{M+r}(x).
\end{aligned}
\ee
Hence
\be{upst5-51}
\begin{aligned}
\sum_{x\in\Z^3}\Xi_M(x)w(x,r)
&\leq \frac{4\alpha}{\kappa}e^{C\alpha/T} \int_0^r d\widetilde r
\sum_{x\in\Z^3} \Xi_M(x)\Xi_{M+\widetilde r}(x)\\
&\leq \frac{4\alpha}{\kappa}e^{C\alpha/T}
\frac{1}{R_\kappa^2}\int_{k R_\kappa}^{(k+1)R_\kappa} ds
\int_{k R_\kappa}^{(k+1)R_\kappa} d\widetilde s \int_0^r d\widetilde r\,
p_{12T\onek+\frac{s+\widetilde s -2u+2M+\widetilde r}\kappa}(0)\\
&\leq \frac{4\alpha}{\kappa}e^{C\alpha/T} \int_0^\infty d\widetilde r\,
p_{\frac{2M+\widetilde r}\kappa}(0)\\
&\leq \frac{4C\alpha}{\sqrt{K}\kappa}e^{C\alpha/T},
\end{aligned}
\ee
where we again use the second inequality of Lemma \ref{approxLaplem}. Substituting
(\ref{upst5-51}) into (\ref{upst5-31}), we arrive at the claim
with $D_{\alpha,T,K}=8\alpha^2Ce^{2C\alpha/T}/\sqrt{K}$.
\epr


\subsubsection{Further reduction of Lemma~\ref{defreez*}}
\label{S5.reddefreez*}

To further estimate the expectation in Lemma \ref{uplem3} from above, we use the
following two lemmas.

\bl{GdHMlem}
Let
\be{GdHM-1}
\Gamma(\beta)
=\limsup_{t\to\infty}\frac1t\log\E_{\,\nu_\rho}\bigg(\exp\bigg[
\beta\int_0^t du\, \big[\xi_u(e_1)-\xi_u(0)\big]^2\bigg]\bigg).
\ee
Then
\be{GdHM-3}
\lim_{\beta\to 0}\frac{\Gamma(\beta)}{\beta}
=2\rho(1-\rho).
\ee
\el

\bpr
The proof is a straightforward adaptation of what is done in G\"artner,
den Hollander and Maillard~\cite{garholmai06}, Lemmas 4.6.8 and 4.6.10.
\epr

\bl{ublem10}
For all $\alpha,T,K,R,\kappa>0$, $u\in[(k-1)R_\kappa,k R_\kappa]$,
$k=1,2,\cdots,\lfloor t/R_\kappa\rfloor$, and all paths $X$,
\be{uplem10-1}
\limsup_{t\to\infty}\frac1{2t} \log\cE^{(5)}_{k,u}(t)
\leq \vartheta_{\alpha,T}\, \rho(1-\rho)\big\|\cK_{k,u}\big\|_1,
\ee
where $\vartheta_{\alpha,T}$ does not depend on $K,R,\kappa,u,k$ or $X$, and
$\vartheta_{\alpha,T}\to 1$ as $T\to\infty$.
\el

\bpr
Using the bound in (\ref{K1bd}) for $\|\cK_{k,u}\|_1$, we find that
\be{}
\kappa\big\|\cK_{k,u}\big\|_1
\leq e^{2C\alpha/T}2\alpha^2\int_0^\infty dr\, p_{12T+2r}(0)
\leq\frac{C\alpha^2 e^{2C\alpha/T}}{\sqrt{T}},
\ee
which tends to zero as $T\to\infty$. Hence, we may apply Lemma~\ref{GdHMlem} to
(\ref{uplem3-1*}) to get the claim.
\epr

At this point we may combine Lemmas~\ref{ublem10} and \ref{st5lem} with (\ref{upst5-11}),
to get
\be{lambkubd}
\lambda_{k,u}^{\kappa,X}
\leq \vartheta_{\alpha,T}\, \rho(1-\rho)\big\|\cK_{k,u}\big\|_1
+\frac{D_{\alpha,T,K}}{2\kappa^2}\rho.
\ee
Note that the upper bound in (\ref{K1bd}) for $\|\cK_{k,u}\|_1$ depends on $X$ only via
its increments on the times interval $[(k-1)R_\kappa,k R_\kappa]$ and that these increments
are i.i.d.\ for $k=1,2,\cdots,\lfloor t/R_\kappa\rfloor$. Hence, combining (\ref{E4-3})
and Lemma \ref{spectralbdlem} with (\ref{lambkubd}) and splitting the resulting expectation
w.r.t.\ $\ES_0^X$ into $\lfloor t/R_\kappa\rfloor$ equal factors with the help of the Markov
property at times $k R_\kappa$, $k=1,2,\cdots,\lfloor t/R_\kappa\rfloor$, we obtain, after
also substituting (\ref{K1bd}),
\be{}
\limsup_{t\to\infty}\frac{\kappa^2}{t}\log \cE^{(4)}_{R,\alpha}(t)
\leq \frac1R \log \cE^{(7)}_{R,\alpha}(\kappa) + \frac{D_{\alpha,T,K}}2 \rho
\ee
with
\be{E7def}
\begin{aligned}
&\cE^{(7)}_{R,\alpha}(\kappa) = \cE^{(7)}_{R,\alpha}(T,K;\kappa)\\
&=\ES_{0}^{X}\Bigg(\exp\Bigg[
\frac{\Theta_{\alpha,T,\rho}}{\kappa^2} \frac{1}{R_\kappa^2}\int_{0}^{R_\kappa} ds
\int_{s}^{R_\kappa} d\widetilde s \int_{-R_\kappa}^0 du \int_0^M dr\,
p_{12T\onek+\frac{s+\widetilde s-2u+2r}{\kappa}}\big(X_{\widetilde s}-X_s\big)\Bigg]\Bigg),
\end{aligned}
\ee
where
\be{Thetadef}
\Theta_{\alpha,T,\rho}
=4\vartheta_{\alpha,T}\alpha^2 e^{2C\alpha/T}\rho(1-\rho)\to 4\alpha^2\rho(1-\rho)
\quad\text{as } T\to\infty.
\ee
Because of (\ref{st5lem-2}), we therefore conclude that the proof of Lemma~\ref{defreez*}
reduces to the following lemma.
\bl{upbdlem12}
For all $\alpha,K>0$,
\be{upbdlem12-1}
\limsup_{\kappa,T,R\to\infty}\frac{1}{R}\log\cE^{(7)}_{R,\alpha}(T,K;\kappa)
\leq \big[6\alpha^2\rho(1-\rho)\big]^2\cP_3.
\ee
\el


\subsection{Small-time cut out, scaling and large deviations}
\label{S5.3}


\subsubsection{Small-time cut out}
\label{S5.3.1}

The proof of Lemma~\ref{upbdlem12} will be reduced to two further lemmas in which
we cut out small times. These lemmas will be proved in Sections
\ref{5.upbdlem12*}--\ref{S5.4}.

For $\epsilon>0$ small, let
\be{mdef}
m=3\epsilon\kappa^3\onek
\ee
and define
\be{E8def}
\begin{aligned}
&\cE^{(8)}_{R,\alpha}(\kappa)=\cE^{(8)}_{R,\alpha}(T,\epsilon;\kappa)\\
&=\ES_{0}^{X}\Bigg(\exp\Bigg[
\frac{\Theta_{\alpha,T,\rho}}{\kappa^2 R_\kappa^2} \int_{0}^{R_\kappa} ds
\int_{s}^{R_\kappa} d\widetilde s \int_{-R_\kappa}^0 du \int_0^m dr\,
p_{12T\onek+\frac{s+\widetilde s-2u+2r}{\kappa}}
\big(X_{\widetilde s}-X_s\big)\Bigg]\Bigg)
\end{aligned}
\ee
and
\be{E9def}
\begin{aligned}
&\cE^{(9)}_{R,\alpha}(\kappa)=\cE^{(9)}_{R,\alpha}(T,\epsilon,K;\kappa)\\
&=\ES_{0}^{X}\Bigg(\exp\Bigg[
\frac{\Theta_{\alpha,T,\rho}}{\kappa^2 R_\kappa^2} \int_{0}^{R_\kappa} ds
\int_{s}^{R_\kappa} d\widetilde s \int_{-R_\kappa}^0 du \int_m^M dr\,
p_{12T\onek+\frac{s+\widetilde s-2u+2r}{\kappa}}
\big(X_{\widetilde s}-X_s\big)\Bigg]\Bigg).
\end{aligned}
\ee
By H\"older's inequality with weights $p,q>1$, $1/p+1/q=1$, we have
\be{}
\cE^{(7)}_{R,\alpha}(\kappa)
= \Big(\cE^{(8)}_{R,\sqrt{q}\alpha}(\kappa)\Big)^{1/q}
\Big(\cE^{(9)}_{R,\sqrt{p}\alpha}(\kappa)\Big)^{1/p}.
\ee
Hence, by choosing $p$ close to $1$, we see that the proof of Lemma \ref{upbdlem12}
reduces to the following lemmas.

\bl{upbdlem12*}
For all $\alpha>0$ and $\epsilon>0$ small enough,
\be{upbdlem12*-1}
\limsup_{\kappa,T,R\to\infty}\frac{1}{R}\log\cE^{(8)}_{R,\alpha}(T,\epsilon;\kappa)
=0.
\ee
\el

\bl{upbdlem13}
For all $\alpha,\epsilon,K>0$ with $0<\epsilon<K$,
\be{upbdlem13-1}
\limsup_{\kappa,T,R\to\infty}\frac{1}{R}\log\cE^{(9)}_{R,\alpha}(T,\epsilon,K;\kappa)
\leq \big[6\alpha^2\rho(1-\rho)\big]^2\cP_3.
\ee
\el

\noindent
Note that in $\cE^{(8)}_{R,\alpha}(\kappa)$ we integrate the transition kernel over
``small'' times $r\in[0,m]$. What Lemma~\ref{upbdlem12*} shows is that the integral
is asymptotically negligible.

\subsubsection{Proof of Lemma~\ref{upbdlem12*}}
\label{5.upbdlem12*}

\bpr
We need only prove the upper bound in (\ref{upbdlem12*-1}). An application of Jensen's
inequality yields
\be{E8eq3}
\begin{aligned}
\cE^{(8)}_{R,\alpha}(\kappa)
\leq\frac{1}{R_\kappa}\int_{0}^{R_\kappa} ds\, \ES_{0}^{X}\Bigg(\exp\Bigg[
\frac{\Theta_{\alpha,T,\rho}}{\kappa^2 R_\kappa}
\int_{0}^{\infty} d\widetilde s \int_{-R_\kappa}^0 du \int_0^m dr\,
p_{12T\onek+2\frac{s-u+r}{\kappa}+\frac{\widetilde s}{\kappa}}
\big(X_{\widetilde s}\big)\Bigg]\Bigg).
\end{aligned}
\ee
Observe that
\be{}
p_{12T\onek+2\frac{s-u+r}{\kappa}+\frac{\widetilde s}{\kappa}}\big(X_{\widetilde s}\big)
=\ES_0^Y\Big(p_{12T\onek+2\frac{s-u+r}{\kappa}}\big(X_{\widetilde s}
+Y_{\widetilde s/\kappa}\big)\Big).
\ee
As in (\ref{Xhat}), let $\widehat X_t=X_t+Y_{t/\kappa}$ and let $\ES_0^{\widehat X}$ denote
expectation w.r.t.\ $\widehat X$ starting at $0$. Then, using Jensen's inequality w.r.t.\
$\ES_0^Y$, we find that
\be{E8eq5}
\begin{aligned}
\cE^{(8)}_{R,\alpha}(\kappa)
\leq\frac{1}{R_\kappa}\int_{0}^{R_\kappa} ds\, \ES_{0}^{\widehat X}\Bigg(\exp\Bigg[
\frac{\Theta_{\alpha,T,\rho}}{\kappa^2 R_\kappa}
\int_{0}^{\infty} d\widetilde s \int_{-R_\kappa}^0 du \int_0^m dr\,
p_{12T\onek+2\frac{s-u+r}{\kappa}}\Big(\widehat X_{\widetilde s}\Big)\Bigg]\Bigg).
\end{aligned}
\ee
For the potential
\be{Vpoten}
V_{s}^\kappa(x)
=\frac{1}{\kappa^2 R_\kappa}\int_{-R_\kappa}^0 du\,\int_0^m dr\,
p_{12T\onek+2\frac{s-u+r}{\kappa}}(x),
\ee
we obtain
\be{}
\begin{aligned}
\Big\|\widehat\cG V_{s}^\kappa\Big\|_\infty
\leq \frac{1}{\kappa^2}\int_0^m dr\, \widehat G_{2T+\frac{r}{3\kappa\onek}}(0)
\leq \frac{3}{\kappa}\onek\int_0^{\epsilon\kappa^2} dr\, \widehat G_{r}(0)
\leq C\sqrt{\epsilon},
\end{aligned}
\ee
where $\widehat \cG$ and $\widehat G$ are the Green operator, respectively, the Green
function corresponding to $\onek\Delta$. Hence, an application of Lemma~\ref{greenoplem}
to (\ref{E8eq5}) yields
\be{E8eq9}
\cE^{(8)}_{R,\alpha}(\kappa)
\leq\big(1-C\Theta_{\alpha,T,\rho}\sqrt{\epsilon}\big)^{-1},
\ee
which, together with (\ref{Thetadef}), leads to the claim for
$0<\epsilon<(4C\rho(1-\rho)\alpha^2)^{-2}$.
\epr

For further comments on Lemma \ref{upbdlem12*}, see the remark
at the end of Section~\ref{S5.4}.

\subsubsection{Scaling, compactification and large deviations}
\label{S5.4}

In this section we prove Lemma~\ref{upbdlem13} with the help of scaling, compactification
and large deviations.

\bpr
Recalling the definition of $m$ in (\ref{mdef}) and $M$ in (\ref{Sregime}), we obtain from
(\ref{E9def}), after appropriate time scaling ($s\to \kappa^2 s$, $\widetilde s\to\kappa^2
\widetilde s$, $u\to\kappa^2 u$ and $r\to 3\kappa^3\onek r$),
\be{E9eq3}
\begin{aligned}
&\cE^{(9)}_{R,\alpha}(\kappa)\\
&=\ES_{0}^{X}\bigg(\exp\bigg[
3\Theta_{\alpha,T,\rho}\onek \frac{1}{R^2}\int_{0}^{R} ds
\int_{s}^{R} d\widetilde s \int_{-R}^0 du \int_{\epsilon}^{K} dr\,
p^{(\kappa)}_{\frac{2T\onek}{\kappa^2}+\frac{s+\widetilde s-2u}{6\kappa}+\onek r}
\Big(X^{(\kappa)}_{s},X^{(\kappa)}_{\widetilde s}\Big)\bigg]\bigg)
\end{aligned}
\ee
with the rescaled transition kernel
\be{}
p_t^{(\kappa)}(x,y)
=\kappa^3 p_{6\kappa^2t}(\kappa x,\kappa y),
\quad x,y\in\Z^3_\kappa=\frac1\kappa \Z^3,
\ee
and the rescaled random walk
\be{}
X_t^{(\kappa)}=\kappa^{-1} X_{\kappa^2 t},
\quad t\in [0,\infty).
\ee
Let $Q$ be a large centered cube in $\R^3$, viewed as a torus, and let $Q^{(\kappa)}
=Q\cap\Z_\kappa^3$. Let $l(Q)$, $l(Q^{(\kappa)})$ denote the side lengths of $Q$ and
$Q^{(\kappa)}$, respectively. Define the periodized objects
\be{upst7-23}
p_t^{(\kappa,Q)}(x,y) =\sum_{k\in\Z^3} p_t^{(\kappa)}\Big(x,y+\frac{k}{\kappa}\,
l\big(Q^{(\kappa)}\big)\Big)
\ee
and
\be{upst7-21}
X_t^{(\kappa,Q)}= X_t^{(\kappa)}
\quad\text{mod}\big(Q^{(\kappa)}\big).
\ee
Clearly,
\be{}
p_t^{(\kappa)}\Big(X_s^{(\kappa)},X_{\widetilde s}^{(\kappa)}\Big)
\leq p_t^{(\kappa,Q)}\Big(X_s^{(\kappa,Q)},X_{\widetilde s}^{(\kappa,Q)}\Big).
\ee
Let $\beta=(\beta_t)_{t\geq 0}$ be Brownian motion on the torus $Q$
with generator $\Delta_{\R^3}$ and transition kernel
\be{upst7-27}
p_t^{(G,Q)}(x,y)
=\sum_{k\in\Z^3}p_t^{(G)}\Big(x,y+k\,l(Q)\Big)
\ee
obtained by periodization of the Gaussian kernel $p_t^{(G)}(x,y)$ defined
in (\ref{Gaussianker}). Fix $\theta>1$ (arbitrarily close to $1$). Then there
exists $\kappa_0=\kappa_0(\theta;\epsilon,K,Q)>0$ such that
\be{}
p_t^{(\kappa,Q)}(x,y) \leq \theta p_t^{(G,Q)}(x,y),
\quad \text{for all } \kappa>\kappa_0
\text{ and } (t,x,y)\in [\epsilon/2,2K]\times Q\times Q.
\ee
Hence, it follows from (\ref{E9eq3}) that there exists
$\kappa_1=\kappa_1(\theta;T,\epsilon,K,R,Q)>0$ such that
\be{E9-bd*}
\cE^{(9)}_{R,\alpha}(\kappa)
\leq \ES_0^X\Bigg(\exp\Bigg[\frac32\theta^2
\Theta_{\alpha,T,\rho}\frac{1}{R}\int_{0}^{R} ds
\int_{0}^{R} d\widetilde s \int_{\epsilon}^{K} dr\,
p^{(G,Q)}_{r}\Big(X^{(\kappa,Q)}_{s},X^{(\kappa,Q)}_{\widetilde s}\Big)\Bigg]\Bigg).
\ee
Applying Donsker's invariance principle and recalling (\ref{Thetadef}),
we find that
\be{}
\begin{aligned}
&\limsup_{\kappa,T\to\infty}\frac1R \log \cE^{(9)}_{R,\alpha}(\kappa)\\
&\quad\leq \frac1R \log \ES_0^\beta\Bigg(\exp\Bigg[6\theta^2
\alpha^2\rho(1-\rho) \frac{1}{R}\int_{0}^{R} ds
\int_{0}^{R} d\widetilde s \int_{\epsilon}^{K} dr\,
p^{(G,Q)}_{r}\big(\beta_{s},\beta_{\widetilde s}\big)\Bigg]\Bigg).
\end{aligned}
\ee
Applying the large deviation principle for the occupation time measures of $\beta$,
we get
\be{upst7-48}
\limsup_{\kappa,T,R\to\infty}\frac1R \log \cE^{(9)}_{R,\alpha}(T,\epsilon;\kappa)
\leq \cP_3^{(Q)}(\theta;\epsilon,K),
\ee
where
\be{upst7-48*}
\cP_3^{(Q)}(\theta;\epsilon,K)
=\sup_{\nu\in\cM_1(Q)}\Bigg[6\theta^2\alpha^2\rho(1-\rho)\int_Q \nu(dx)\int_Q\nu(dy)
\int_{\epsilon}^{K} dr\, p_r^{(G,Q)}(x,y)-S^Q(\nu)\Bigg]
\ee
with large deviation rate function $S^Q\colon\cM_1(Q)\to[0,\infty]$ defined by
\be{upst7-49}
S^Q(\mu) =
\begin{cases}\|\nabla_{\R^3}f\|_{2}^2
&\text{if }
\mu\ll dx\text{ and } \sqrt{\frac{d\mu}{dx}}=f(x)\text{ with }f\in H_{\per}^1(Q),\\
\infty & \text{otherwise,}
\end{cases}
\ee
where $\cM_1(Q)$ is the space of probability measures on $Q$, and $H_{\per}^1(Q)$ denotes
the space of functions in $H^1(Q)$ with periodic boundary conditions. By \cite{garhol04},
Lemma 7.4, we have
\be{upst7-81}
\limsup_{Q\ua\R^3}\cP^{(Q)}_3(\theta;\epsilon,K)
\leq \cP_3(\theta;\epsilon,K)
\ee
with
\be{upst7-83}
\begin{aligned}
&\cP_3(\theta;\epsilon,K)\\
&=\sup_{{f\in H^1(\R^3)} \atop {\|f\|_2=1}}\bigg[6\,\theta^2\alpha^2\rho(1-\rho)\,
\int_{\R^3} dx\, f^2(x) \int_{\R^3} dy\, f^2(y)\int_{\epsilon}^{K} dr\,
p_r^{(G)}(x,y)-\big\|\nabla_{\R^3}f\big\|_{L^2(\R^3)}^{2}\bigg]\\
&\leq\sup_{{f\in H^1(\R^3)} \atop {\|f\|_2=1}}
\bigg[6\,\theta^2\alpha^2\rho(1-\rho)\int_{\R^3} dx\, f^2(x) \int_{\R^3} dy\, f^2(y) \int_{0}^{\infty} dr\,
p_r^{(G)}(x,y)-\big\|\nabla_{\R^3}f\big\|_{L^2(\R^3)}^{2}\bigg]\\
&=\big[6\,\theta^2\alpha^2\rho(1-\rho)\big]^2 \cP_3.
\end{aligned}
\ee
Combining (\ref{upst7-48}) and (\ref{upst7-83}), and letting $\theta\da 1$,
we  arrive at the claim of Lemma \ref{upbdlem13}.
\epr

This, after a long struggle by the authors and considerable patience on the side
of the reader, completes the proof of the upper bound in Proposition \ref{mainLem1*}.

\noindent{\it Remark.} The reader might be surprised that the expression in the l.h.s.\
of (\ref{upbdlem12*-1}) does not only vanish in the limit as $\epsilon\da 0$ but
vanishes for \emph{all} $\epsilon>0$ sufficiently small. This fact is closely related
to the observation that
\be{}
\cP_3\big(\pi^3\big)=0
\quad\text{whereas}\quad
\cP_3(\infty)=\cP_3>0
\ee
with
\be{}
\cP_3(\epsilon)
=\sup_{{f\in H^1(\R^3)} \atop {\|f\|_2=1}}
\Bigg[\int_{\R^3} dx\, f^2(x)\int_{\R^3} dy\, f^2(y)\int_0^\epsilon dr\,
p_r^{(G)}(x-y) - \big\|\nabla_{\R^3} f\big\|_2^2\Bigg].
\ee
Indeed, given a potential $V\geq 0$ with $\|\cG_{\R^3}V\|_\infty <1/2$,
where $\cG_{\R^3}$ denotes the Green operator associated with $\Delta_{\R^3}$,
the method used in the proof of Lemma \ref{upbdlem12*} leads to
\be{}
\lim_{R\to\infty}\frac1R\log \ES_0^\beta \Bigg(\exp\Bigg[
\frac1R\int_0^R ds \int_0^R d\widetilde s\,\,
V(\beta_{\widetilde s}-\beta_{s})\Bigg]\Bigg)
=0.
\ee
On the other hand, the large deviation principle for the occupation time measures of
$\beta$ shows that this limit coincides with
\be{}
\sup_{{f\in H^1(\R^3)} \atop {\|f\|_2=1}}
\Bigg[\int_{\R^3} dx\, f^2(x)\int_{\R^3} dy\, f^2(y)\, V(x-y)
- \big\|\nabla_{\R^3} f\big\|_2^2\Bigg].
\ee
But, for $0<\epsilon<\pi^3$ the potential
\be{}
V_\epsilon(x)=\int_0^\epsilon dr\, p_r^{(G)}(x)
\ee
satisfies the assumption $\|\cG_{\R^3}V_\epsilon\|_\infty< 1/2$, implying
$\cP_3(\pi^3)=0$.


\section{Higher moments}
\label{S6}

In this last section we explain how to extend the proof of Theorem~\ref{main} to
higher moments $p\geq 2$. Sections \ref{S6.1}--\ref{S6.3} parallel Sections \ref{S2.1},
\ref{S3.2}, \ref{S4} and \ref{S5}.

\subsection{Two key propositions}
\label{S6.1}

Our starting point is the Feynman-Kac representation for
the $p$-th moment,
\be{}
\E_{\,\nu_\rho}\Big(u(0,t)^p\Big)
=\E_{\,\nu_\rho;0}^{(p)}\Bigg(\exp\Bigg[\int_0^t ds\, \sum_{j=1}^{p}
\xi_s\big(X_{\kappa s}^j\big)\Bigg]\Bigg),
\ee
where $X^1,\cdots,X^p$ are independent simple random walks on $\Z^3$ starting at
$0$ and $\E_{\,\nu_\rho;x}^{(p)}$ denotes expectation w.r.t.\
$\P_{\,\nu_\rho;x}^{(p)}=\P_{\,\nu_\rho}\otimes\PS_{x_1}^{X^1}\otimes\cdots\otimes\PS_{x_p}^{X^p}$,
$x=(x_1,\cdots,x_p)\in(\Z^3)^p$.

The arguments in Sections \ref{S2} and \ref{S3} easily extend to this more general
case by replacing $Z$, $\cA$, $(\cS_t)_{t\geq 0}$, $\phi$ and $\psi$ by their $p$-dimensional
analogues $Z^{(p)}$, $\cA^{(p)}$, $(\cS_t^{(p)})_{t\geq 0}$, $\phi^{(p)}$ and $\psi^{(p)}$.
To be precise, consider the Markov process
\be{}
Z_t^{(p)}=\big(\xi_{t/\kappa},X_t^1,\cdots,X_t^{p}\big)
\quad\text{on } \Omega\times(\Z^3)^p
\ee
with generator
\be{}
\cA^{(p)}=\frac1\kappa L+\sum_{j=1}^{p} \Delta_j,
\ee
where the lattice Laplacian $\Delta_j$ acts on the $j$-th spatial
variable. Denote by $(\cS^{(p)}_t)_{t\geq 0}$ the associated semigroup. We define
\be{}
\phi^{(p)}(\eta;x_1,\cdots,x_p)
=\sum_{j=1}^{p}\phi(\eta,x_j)
=\sum_{j=1}^{p}(\eta(x_j)-\rho)
\ee
and
\be{}
\psi^{(p)}=\int_0^T ds\, \cS_s^{(p)}\phi^{(p)}.
\ee
Then
\be{}
\psi^{(p)}(\eta;x_1,\cdots,x_p)=\sum_{j=1}^{p}\psi(\eta,x_j).
\ee
In this way the proof of Theorem~\ref{main} for $p\geq 2$ reduces to the proof of the
following extension of Propositions~\ref{mainLem2*} and \ref{mainLem1*}.

\bp{pextprop1}
For all $p\in\N$ and $\alpha\in\R$,
\be{}
\begin{aligned}
&\limsup_{t,\kappa,T\to\infty}\frac{\kappa^2}{pt}
\log\E_{\,\nu_\rho;0}^{(p)}\Bigg(\exp\Bigg[\alpha\,\int_0^t ds\,
\bigg[\Big(e^{-\frac{1}{\kappa} \psi^{(p)}} \cA^{(p)} e^{\frac{1}{\kappa}\psi^{(p)}}\Big)
-\cA^{(p)}\Big(\frac{1}{\kappa}\psi^{(p)}\Big)\bigg]\big(Z_s^{(p)}\big)\Bigg]\Bigg)\\
&\quad\leq\frac{\alpha}{6}\,\rho(1-\rho)G.
\end{aligned}
\ee

\ep
\bp{pextprop2}
For all $p\in\N$ and $\alpha>0$,
\be{}
\lim_{t,\kappa,T\to\infty}
\frac{\kappa^2}{pt}\log \E_{\,\nu_\rho;0}^{(p)}\Bigg(\exp\Bigg[
\frac{\alpha}{\kappa}\int_0^t ds\,
\Big(\cS_{T}^{(p)}\phi^{(p)}\Big)\big(Z_s^{(p)}\big)\Bigg]\Bigg)
= \big[6\alpha^2\rho(1-\rho)p\big]^2\cP_3.
\ee
\ep

\noindent
Proposition~\ref{pextprop1} has already been proven for all $p\in\N$ in \cite{garholmai06},
Proposition 4.4.2 and Section 4.8.


\subsection{Lower bound in Proposition \ref{pextprop2}}
\label{S6.2}

We use the following analogue of the variational representation (\ref{lwst1-3})
\be{}
\begin{aligned}
&\lim_{t\to\infty}\frac{1}{t}\log\E_{\,\nu_\rho;0}^{(p)}\Bigg(\exp\Bigg[
\frac\alpha\kappa\int_0^t ds\, \Big(\cS_{T}^{(p)}\phi^{(p)}\Big)
\big(Z_s^{(p)}\big)\Bigg]\Bigg)\\
&\qquad=\sup_{{F^{(p)}\in\cD(\cA^{(p)})}\atop{\|F^{(p)}\|_{L^2(\mu_\rho^{p})}=1}}
\iint_{\Omega\times(\Z^3)^{p}}d\mu_\rho^{p}
\bigg[\frac{\alpha}{\kappa} \Big(\cS_{T}^{(p)}\phi^{(p)}\Big)
\big(F^{(p)}\big)^2 +F^{(p)} \cA^{(p)} F^{(p)}\bigg].
\end{aligned}
\ee
To obtain the appropriate lower bound, we use test functions $F^{(p)}$ of the form
\be{}
F^{(p)}(\eta;x_1,\cdots,x_p)
=F_1(\eta)F_2(x_1)\cdots F_2(x_p),
\ee
where $F_1$, $F_2$ and $F=F^{(1)}$ are the same as in (\ref{F1def}), (\ref{F2def}) and
(\ref{testF}), respectively. One easily checks that
\be{}
\begin{aligned}
&\frac{\kappa^2}{p}\iint_{\Omega\times(\Z^3)^{p}}d\mu_\rho^{p}
\bigg[\frac{\alpha}{\kappa} \Big(\cS_{T}^{(p)}\phi^{(p)}\Big)
\big(F^{(p)}\big)^2 +F^{(p)} \cA^{(p)} F^{(p)}\bigg]\\
&\quad=\frac{(p\kappa)^2}{p^2}\iint_{\Omega\times \Z^3}d\mu_\rho
\bigg[\frac{p\alpha}{p\kappa} \Big(\cS_{T}\phi\Big)
F^2 +F \bigg(\frac{1}{p\kappa} L+ \Delta\bigg) F \bigg].
\end{aligned}
\ee
But this is $1/p^2$ times the expression in Lemma \ref{testfunct} with $\alpha$ and $\kappa$
replaced by $p\alpha$ and $p\kappa$, respectively. Hence, we may use the lower bounds for this
expression in Section~\ref{S4} to arrive at the lower bound in Proposition~\ref{pextprop2}.


\subsection{Upper bound in Proposition \ref{pextprop2}}
\label{S6.3}

1. Freezing and defreezing can be done in the same way as in Section
\ref{S5.1}, but with $V(\eta,x)$ in (\ref{upst2-3}) replaced by
\be{}
V^{(p)}(\eta,x)
=\frac\alpha\kappa\sum_{y\in\Z^3}\Bigg(\sum_{j=1}^p p_{6T\onek}\big(x_j,y\big)\Bigg)
\big(\eta(y)-\rho\big).
\ee
This leads to the analogues of Lemmas \ref{freez} and \ref{defreez} along
the lines of Sections \ref{S5.freez} and \ref{S5.defreez}.

\noindent
2. Considering
\be{}
V_{k,u}^{(p)}(\eta)
=\frac1{R_\kappa}\int_{k R_\kappa}^{(k+1)R_\kappa}ds \sum_{y\in \Z^3}
\Bigg(\sum_{j=1}^{p}p_{6T\onek+\frac{s-u}{\kappa}}\big(X_s^j,y\big)\Bigg)
\big(\eta(y)-\rho\big)
\ee
and
\be{}
\cE_{R,\alpha}^{(4,p)}(t)
=\E_{\,\nu_\rho;0}\Bigg(\exp\Bigg[\frac{\alpha}{\kappa}\sum_{k=1}^{\lfloor t/R_\kappa\rfloor}
\int_{(k-1)R_\kappa}^{k R_\kappa} du\, V_{k,u}^{(p)}\big(\xi_{u/\kappa}\big)\Bigg]\Bigg)
\ee
instead of (\ref{Vpot}--\ref{E4-3}), the proof reduces to the following analogue of
Lemma~\ref{defreez*}.

\bl{limfin}
For each $\alpha>0$,
\be{}
\limsup_{R\to\infty}\limsup_{t,\kappa,T\to\infty}\frac{\kappa^2}{pt}
\log\cE_{R,\alpha}^{(4,p)}(t)
\leq \big[6\alpha^2\rho(1-\rho)p\big]^2\cP_3.
\ee
\el

\noindent
3. The proof of Lemma \ref{limfin} follows the lines of Sections \ref{S5.defreez*}--\ref{S5.3}.
The spectral bound is essentially the same as in Section \ref{S5.2}.
In Lemma \ref{spectralbdlem} we have to replace $V_{k,u}$ by $V_{k,u}^{(p)}$
and $\lambda_{k,u}$ by
\be{}
\lambda_{k,u}^{(p)}=\lim_{t\to\infty}\frac1t\log\E_{\,\nu_\rho}\Bigg(\exp\Bigg[
\frac\alpha\kappa\int_0^t ds\, V_{k,u}^{(p)}\big(\xi_{s/\kappa}\big)\Bigg]\Bigg).
\ee
Subsequently, we replace $V_{k,u}$ by $V_{k,u}^{(p)}$ in (\ref{upst5-3}), to obtain functions
$\widehat \phi^{(p)}$, $\widehat \psi^{(p)}$ replacing (\ref{upst5-3}--\ref{psihatdef}), and
\be{}
\Xi_r^{(p)}(x)=\frac{1}{R_\kappa}\int_{k R_\kappa}^{(k+1)R_\kappa} ds\,
\sum_{j=1}^{p} p_{6T\onek+\frac{s-u+r}{\kappa}}\big(X_s^j,x\big)
\ee
replacing (\ref{uplem3-3}).
Then, in the analogue of Lemma \ref{uplem3}, instead of (\ref{K1bd}) we get the bound
\be{}
\big\|\cK_{k,u}^{(p)}\big\|_1
\leq e^{2C\alpha/T}\frac{2\alpha^2}{\kappa^2 R_\kappa^2}
\int_{k R_\kappa}^{(k+1)R_\kappa} ds\int_{k R_\kappa}^{(k+1)R_\kappa} d\widetilde s
\int_0^M dr \sum_{i,j=1}^{p}
p_{12T\onek+\frac{s+\widetilde s-2u+2r}{\kappa}}\big(X_{\widetilde s}^j-X_s^i\big)
\ee
along the line of Section \ref{S5.uplem}. Similarly, the proof of the analogue of
Lemma \ref{st5lem} follows the argument in Section \ref{S5.st5lem}, leading to a
reduction of Lemma \ref{limfin} to the analogue of Lemma \ref{upbdlem12}, as in
Section \ref{S5.reddefreez*}.

\noindent
4. To make the small-time cut-off, instead of (\ref{E8def}) we consider
\be{E8p}
\begin{aligned}
&\cE_{R,\alpha}^{(8,p)}(\kappa)\\
&=\ES_0^X\Bigg(\exp\Bigg[\frac{\Theta_{\alpha,T,\rho}}{\kappa^2 R_\kappa^2}
\int_{0}^{R_\kappa}ds\int_{s}^{R_\kappa}d\widetilde s\int_{-R_\kappa}^{0}du\int_0^m dr\,
\sum_{i,j=1}^{p} p_{12T\onek+\frac{s+\widetilde s-2u+2r}{\kappa}}
\big(X_{\widetilde s}^j-X_s^i\big)
\Bigg]\Bigg).
\end{aligned}
\ee
Using the Chapman-Kolmogorov equation, we see that
\be{}
\begin{aligned}
&\int_{0}^{R_\kappa} ds\int_{0}^{R_\kappa} d\widetilde s\sum_{i,j=1}^p
p_{12T\onek+\frac{s+\widetilde s-2u+2r}{\kappa}}\big(X_{\widetilde s}^j-X_s^i\big)\\
&\qquad=\sum_{z\in\Z^3}\Bigg(\sum_{i=1}^{p}\int_{0}^{R_\kappa} ds\,
p_{6T\onek+\frac{s-u+r}{\kappa}}\big(X_s^i,z\big)\Bigg)^2\\
&\qquad\leq p\sum_{z\in\Z^3}\sum_{i=1}^{p}\Bigg(\int_{0}^{R_\kappa} ds\,
p_{6T\onek+\frac{s-u+r}{\kappa}}\big(X_s^i,z\big)\Bigg)^2\\
&\qquad=p\sum_{i=1}^{p}\int_0^{R_\kappa}ds\int_0^{R_\kappa}d\widetilde s\,
p_{12T\onek+\frac{s+\widetilde s-2u+2r}{\kappa}}\big(X_{\widetilde s}^i-X_s^i\big).
\end{aligned}
\ee
Substituting this into the r.h.s.\ of (\ref{E8p}) and applying H\"older's inequality
for the $p$ exponential factors, we reduce the problem to the consideration of
a single random walk and can proceed as in Section \ref{S5.3}, leading to the analogues
of Lemmas \ref{upbdlem12*}--\ref{upbdlem13}.

\noindent
5. The proof of the analogue of Lemma \ref{upbdlem12*} runs along the line of
Section \ref{5.upbdlem12*}. To prove the analogue of Lemma \ref{upbdlem13}, instead
of (\ref{E9def}) we consider
\be{E9p}
\begin{aligned}
&\cE^{(9,p)}_{R,\alpha}(\kappa)\\
&=\ES_0^X\Bigg(\exp\Bigg[\frac{\Theta_{\alpha,T,\rho}}{\kappa^2 R_\kappa^2}
\int_{0}^{R_\kappa}ds\int_{s}^{R_\kappa}d\widetilde s\int_{-R_\kappa}^{0}du\int_m^M dr\,
\sum_{i,j=1}^{p} p_{12T\onek+\frac{s+\widetilde s-2u+2r}{\kappa}}
\big(X_{\widetilde s}^j-X_s^i\big)
\Bigg]\Bigg).
\end{aligned}
\ee
As in Section \ref{S5.4}, this leads to
\be{E9pbd}
\begin{aligned}
&\limsup_{\kappa,T,R\to\infty}\frac1{pR} \log\cE^{(9,p)}_{R,\alpha}(\kappa)\\
&\quad\leq\frac1p\sup_{{\nu_i\in\cM_1(Q)} \atop {1\leq i\leq p}}
\Bigg[6\theta^2\alpha^2\rho(1-\rho)\sum_{i,j=1}^p\int_Q\nu_i(dx)\int_Q\nu_j(dy)
\int_\epsilon^Kdr\,p_r^{(G,Q)}(x,y)-\sum_{i=1}^{p}S^Q(\nu_i)\Bigg]
\end{aligned}
\ee
instead of (\ref{upst7-48}--\ref{upst7-48*}). Now we can proceed similarly as in
\cite{garhol04}, Lemma 7.3. Consider the Fourier transforms $\widehat\nu_i$ of
the measures $\nu_i\in\cM_1(Q)$ defined by
\be{}
\widehat\nu_j(k)
=\int_Q e^{-i(2\pi/l(Q)) k\cdot x} \nu_j(dx),
\quad k\in\Z^3, j=1,\cdots,p.
\ee
The transition kernel $p^{(G,Q)}$ admits the Fourier representation
\be{}
p_{r}^{(G,Q)}(x) = \frac{1}{l(Q)^3}\sum_{k\in\Z^3}
e^{-(2\pi/l(Q))^2 |k|^2r}e^{-i(2\pi/l(Q)) k\cdot x},
\quad (x,t)\in Q\times (0,\infty).
\ee
Therefore we may write
\be{}
\int_Q\nu_i(dx)\int_Q\nu_j(dy)\, p_r^{(G,Q)}(x,y)
=\frac{1}{l(Q)^3}\sum_{k\in\Z^3}
e^{-(2\pi/l(Q))^2 |k|^2r}\widehat \nu_i(k) \overline{\widehat\nu_j(k)}.
\ee
Using that
\be{}
\mbox{Re}\Big(\widehat \nu_i(k) \overline{\widehat\nu_j(k)}\Big)
\leq\frac12 \big|\widehat \nu_i(k)\big|^2+ \frac12  \big|\widehat \nu_j(k)\big|^2,
\ee
we obtain
\be{}
\begin{aligned}
&\int_Q\nu_i(dx)\int_Q\nu_j(dy)\, p_r^{(G,Q)}(x,y)\\
&\quad\leq \frac12 \int_Q\nu_i(dx)\int_Q\nu_i(dy)\, p_r^{(G,Q)}(x,y)
+\frac12 \int_Q\nu_j(dx)\int_Q\nu_j(dy)\, p_r^{(G,Q)}(x,y).
\end{aligned}
\ee
Therefore the term inside the square brackets in the r.h.s.\ of (\ref{E9pbd})
does not exceed
\be{}
\sum_{i=1}^p \Bigg[6\theta^2\alpha^2\rho(1-\rho)p\int_Q \nu_i(dx) \int_Q \nu_i(dy)
\int_\epsilon^K dr\, p_r^{(G,Q)}(x,y)-S^Q(\nu_i)\Bigg],
\ee
and we arrive at
\be{E9pbd*}
\limsup_{\kappa,T,R\to\infty}\frac{p}R \log\cE^{(9,p)}_{R,\alpha}(\kappa)
\leq \cP_3^{(Q,p)}(\theta;\epsilon,K)
\ee
with
\be{}
\cP_3^{(Q,p)}(\theta;\epsilon,K)
=\sup_{\nu\in\cM_1(Q)}
\Bigg[6\theta^2\alpha^2\rho(1-\rho)p\int_Q\nu(dx)\int_Q\nu(dy)
\int_\epsilon^Kdr\,p_r^{(G,Q)}(x,y)-S^Q(\nu)\Bigg],
\ee
which is the analogue of (\ref{upst7-48}--\ref{upst7-48*}) for $p\geq 2$.
The rest of the proof can be easily obtained from the analogues of
(\ref{upst7-81}--\ref{upst7-83}).



\begin{thebibliography}{99}

\bibitem{donvar83}
M.D.\ Donsker and S.R.S.\ Varadhan,
Asymptotics for the polaron,
Comm.\ Pure Appl.\ Math.\ 36 (1983) 505--528.

\bibitem{garhol04}
J.\ G\"artner and F.\ den Hollander,
Intermittency in a catalytic random medium,
Ann.\ Probab.\ 34 (2006) 2219--2287.

\bibitem{garholmai06}
J.\ G\"artner, F.\ den Hollander and G.\ Maillard,
Intermittency on catalysts: symmetric exclusion,
Electronic J.\ Probab.\ 12 (2007) 516--573.

\bibitem{garholmaiHvW}
J.\ G\"artner, F.\ den Hollander and G.\ Maillard,
Intermittency on catalysts,
in: {\it Trends in Stochastic Analysis} (eds.\ J.\ Blath, P.\ M\"orters and
M.\ Scheutzow),
London Mathematical Society Lecture Note Series 353, Cambridge University Press,
Cambridge, 2009.

\bibitem{garkon04}
J.\ G\"artner and W.\ K\"onig,
The parabolic Anderson model,
in: \emph{Interacting Stochastic Systems} (eds.\ J.-D.\ Deuschel and A.\ Greven),
Springer, Berlin, 2005, pp.\ 153--179.

\bibitem{hol00}
F.\ den Hollander,
\emph{Large Deviations},
Fields Institute Monographs 14, American Mathematical Society, Providence, RI, 2000.

\bibitem{lig85}
T.M.\ Liggett,
\emph{Interacting Particle Systems},
Grundlehren der Mathematischen Wissenschaften 276, Springer, New York, 1985.

\end{thebibliography}
\end{document}